\documentclass{amsart}
\usepackage{mathrsfs}
\usepackage{amsfonts}
\usepackage[all]{xy}
\usepackage{latexsym,amsfonts,amssymb,amsmath,amsthm,amscd,amsmath,graphicx}
\usepackage[all]{xy}
\usepackage{ifthen}
\usepackage[T1]{fontenc}
\usepackage{subsupscripts}

\newboolean{workmode}
\setboolean{workmode}{false}
\newboolean{sections}
\setboolean{sections}{true}

\input xy
\xyoption{all}

\makeatletter
\def\th@plain{%
  \thm@notefont{}% same as heading font
  \itshape % body font
}
\def\th@definition{%
  \thm@notefont{}% same as heading font
  \normalfont % body font
}
\makeatother

\newcommand{\hypref}[2]{{\hyperref[#1]{#2~\ref{#1}}}}

\newcommand{\ifwork}[1]{\ifthenelse{\boolean{workmode}}{#1}{}}
\newcommand{\comment}[1]{}
\newcommand{\mute}[1]{}
\newcommand{\printname}[1]{}

\ifwork{
\renewcommand{\comment}[1]{{\marginpar{*}\ \scriptsize{#1}\ }}
\renewcommand{\mute}[1]{{\scriptsize \ #1\ }\marginpar{\scriptsize muted}}
\renewcommand{\printname}[1]
    {\smash{\makebox[0pt]{\hspace{-1.0in}\raisebox{8pt}{\tiny #1}}}}
}

\newcommand{\ifsection}[2]{\ifthenelse{\boolean{sections}}{#1}{#2}}

\theoremstyle{plain}
\ifsection{
    \newtheorem{theorem}{Theorem}[section]
}
{
    
}

\theoremstyle{definition}

\newtheorem{thm}{Theorem}[section]
\newtheorem{cor}[thm]{Corollary}
\newtheorem{lem}[thm]{Lemma}
\newtheorem{prop}[thm]{Proposition}
\theoremstyle{definition}
\newtheorem{defn}[thm]{Definition}

\theoremstyle{remark}
\newtheorem{rem}[thm]{Remark}
\numberwithin{equation}{section}
\newcommand{\sA}{{\mathcal A}}

\newcommand{\sE}{{\mathcal E}}
\newcommand{\sF}{{\mathcal F}}
\newcommand{\sG}{{\mathcal G}}
\newcommand{\sH}{{\mathcal H}}
\newcommand{\sI}{{\mathcal I}}

\newcommand{\sK}{{\mathcal K}}

\newcommand{\sM}{{\mathcal M}}

\newcommand{\sO}{{\mathcal O}}
\newcommand{\sP}{{\mathcal P}}

\newcommand{\sV}{{\mathcal V}}

\newcommand{\A}{{\mathbb A}}

\newcommand{\F}{{\mathbb F}}

\newcommand{\N}{{\mathbb N}}
\renewcommand{\P}{{\mathbb P}}

\newcommand{\Z}{{\mathbb Z}}

\newcommand{\x}{\xrightarrow}

\numberwithin{equation}{section}

\begin {document}
\topmargin= -.2in \baselineskip=20pt

\title{Equivariant special $L$-values of abelian $t$-modules}
\author{Jiangxue Fang}
\address{Department of Mathematics,
Capital Normal University, Beijing 100148, P.R. China} \email{fangjiangxue@gmail.com}
\date{}

\maketitle
%\date{*}%
%\dedicatory{}%
%\commby{}%
% ---------------------------
%\section{questions}
%1. Let $G$ be a finite group. Let $A=\F_q[t]$. For any finite $A[G]$-module $M$, does there exist an element $g\in A[G]$ such that
%for any character $\chi:G\to F^\times$ from $G$ to a filed $F$, we have
%$$\chi(g)=|M\otimes_{A[G]}F[t]|_{F[t]}.$$

%2. Let $L/K$ be a Galois extension of global fields and let $\sO_L$ and $\sO_K$ be the rings of integers of $L$ and $K$, respectively.
%Let $G={\rm Gal}(L/K)$. Under what conditions, $\sO_L$ is a locally free $\sO_K[G]$-module. For which $\mathfrak p\in{\rm Spec}\sO_L$,
%$\sO_L/\mathfrak p\sO_L$ is a locally free $\sO_K/\mathfrak p[G]$-module?

\begin{abstract}
We prove a formula of the equivariant $\infty$-adic special $L$-values of abelian $t$-modules. This gives function field analogues of the
equivariant class number formula. As an application, we calculate the special values
of Artin $L$-functions for Galois representations.
\end{abstract}

\section{{Introduction and statement of the main results}}
\noindent Let $k$ be a finite field of $q$-elements and of characteristic $p$.
%Let $A$ be a $k$-algebra. For any $k$-scheme $X$, let $X\otimes A=X\times_k{\rm Spec}\;A$. Let $G$ be a finite abelian group of order prime to $q$.
%For any $\sO_X$-module $\sF$, denote by $\sF\otimes A$ the pull back of $\sF$ via the projection $X\otimes A\to X$.
Let $F$ be a finite field extension of $k$.
For any $k$-module $M$ and any variable $t$, let $M[t]=M\otimes_kk[t]$, $M[[t^{-1}]]=M\widehat\otimes_kk[[t^{-1}]]$ and $M((t^{-1}))=M\widehat\otimes_kk((t^{-1}))$.
\begin{defn}
(1) For any finite $F[t]$-module $M$, the characteristic polynomial $|M|_{F[t]}$ of $M$ is defined to be
$$|M|_{F[t]}={\rm det}_{F[t]}\Big(1\otimes t-t\otimes 1,\;M\otimes_FF[t]\Big),$$
where the $F[t]$-module structure on $M\otimes_FF[t]$ is given by $F[t]$.

(2) Let $H$ be a finite abelian group of order prime to $p$ and $M$ a finite $k[t][H]$-module.
Then $k[H]\simeq\prod_{i=1}^rk_i$ for some finite fields $k_i$ and $M\simeq\oplus_{i=1}^rM_i$ for some finite $k_i[t]$-module $M_i$.
Define
$$|M|_{k[t][H]}=\prod_{i=1}^r|M_i|_{k_i[t]}\in\prod_{i=1}^rk_i[t]=k[t][H].$$

(3) Let $H$ be a finite abelian group and $M$ a finite $k[t][H]$-module which is free over $k[H]$. Define
$$|M|_{k[t][H]}={\rm det}_{k[t][H]}\Big(1\otimes t-t\otimes 1,\;M\otimes_kk[t]\Big),$$
where the $k[t]$-module structure on $M\otimes_kk[t]$ is given by $k[t]$ and the $k[H]$-module structure is given by $M$.

\end{defn}

Let $K$ be a finite field extension of $k(t)$ and $L$ a finite Galois extension of $K$ with Galois group $G=G(L/K)$.
Let $\sO_K$ and $\sO_L$ be the integral closures of $k[t]$ in $K$ and $L$, respectively. Let $z$ be the image of $t$ in $\sO_K$ and $\sO_L$. For any maximal ideal $\mathfrak p$ of $\sO_K$, choose a prime ideal $\mathfrak P$ of $\sO_L$ above $\mathfrak p$. Let $G_\mathfrak P$ be the decomposition group and $I_\mathfrak P$ the inertia group of $\mathfrak P$ over $\mathfrak p$. Then we have an isomorphism $G_\mathfrak P/I_\mathfrak P\simeq G(\kappa_\mathfrak P/\kappa_\mathfrak p)$ onto the Galois group of the residue field extension $\kappa_\mathfrak P/\kappa_\mathfrak p$.
Then $G_\mathfrak P/I_\mathfrak P$ is generated by the Frobenius element ${\rm Frob}_\mathfrak P$ whose image in $G(\kappa_\mathfrak P/\kappa_\mathfrak p)$
is the $|\kappa_\mathfrak p|$-th power map. We have
$$|\kappa_\mathfrak p|_{k[t]}=\prod_{\delta\in G(\kappa_\mathfrak p/k)}(t-\delta(z))\in k[t]\subset k((t^{-1})).$$
Thus $\frac{1}{|\kappa_\mathfrak p|_{k[t]}}\in t^{-1}k[[t^{-1}]].$
Let $\chi:G\to F^\times$ be a character. Define the $\infty$-adic $L$-value of $\chi$ at $1$ by the convergent infinite product
$$L(1,\,\chi)=\prod_{\mathfrak p\in{\rm Max}(\sO_K)}\Big(1-\frac{\chi(\mathfrak p)}{|\kappa_\mathfrak p|_{k[t]}}\Big)^{-1}\in 1+t^{-1}F[[t^{-1}]],$$
where ${\rm Max}(\sO_K)$ is the set of maximal ideals of $\sO_K$ and where we define $\chi(\mathfrak p)$ as follows:
 $$\chi(\mathfrak p)=\left\{\begin{array}{ll}\chi({\rm Frob}_\mathfrak P)&\hbox{ if }\chi|_{I_\mathfrak P}=1,\\
 0&\hbox{ if }\chi|_{I_\mathfrak P}\neq1.\end{array}\right.$$
Generally, let $\rho:G\to{\rm GL}_F(V)$ be a representation of $G$ on a finite dimensional $F$-vector space $V$.
%\begin{rem}
%Since $\rho$ factor through a $F_0$-representation of $G$ for some finite subextension $F_0$ of $k$ in $F$. From now on, we assume $F$ is a finite extension of $k$ for simplicity.
%\end{rem}
\begin{defn}\label{3161}
For any positive integer $n$, define the $\infty$-adic $L$-value of $\rho$ at $n$ by the convergent infinite product
$$L(n,\,\rho)=\prod_{\mathfrak p\in{\rm Max}(\sO_K)}{\rm det}_{F[[t^{-1}]]}\Big(1-\frac{\rho({\rm Frob}_\mathfrak P)}{|\kappa_\mathfrak p|^n_{k[t]}},\;V^{I_\mathfrak P}[[t^{-1}]]\Big)^{-1}\in 1+t^{-1}F[[t^{-1}]].$$
\end{defn}

\begin{defn}\label{3151}
The $n$-th tensor power of Carlitz module is the functor
$$C^{\otimes n}:\{\sO_K\hbox{-algebras}\}\to\{k[t]\hbox{-modules}\}$$
that associates each $\sO_K$-algebra $B$ to a $k[t]$-module $C^{\otimes n}(B)$ whose underlying $k$-vector space is $B^n$ and whose
$k[t]$-module structure is given by
$$C^{\otimes n}:k[t]\to{\rm End}_k(B^n),\;\;C^{\otimes n}(t)(x_1,\ldots,x_{n-1},x_n)=(zx_1+x_2,\ldots,zx_{n-1}+x_n,zx_n+x_1^q)$$
for any $x_1,\ldots,x_n\in B.$
\end{defn}
Let $V^*={\rm Hom}_F(V,\,F)$. Define a right $F[G]$-module structure on $V^*$ by
$$(\varphi g)(v)=\varphi(\rho(g)(v))\hbox{ for any }\,\varphi\in V^*,\,g\in G\hbox{ and }v\in V.$$ For any left $k[t][G]$-module $M$, we get a $F[t]$-module $V^*\otimes_{k[G]}M$. Thus we get two finite $F[t]$-modules $V^*\otimes_{k[G]}C^{\otimes n}(\sO_L/\mathfrak p\sO_L)$ and $V^*\otimes_{k[G]}{\rm Lie}(C^{\otimes n})(\sO_L/\mathfrak p\sO_L)$.
Let $_{I_\mathfrak P}V$ be the sub-$F[G]$-module of $V$ generated by $gv-v$ for any $g\in I_\mathfrak P$ and $v\in V$ and let $V_{I_\mathfrak P}=V/_{I_\mathfrak P}V$.
In this paper, we will prove the following lemma.
\begin{lem}\label{3164}
We have
\begin{eqnarray*}
L(n,\,\rho)=\prod_{\mathfrak p\in{\rm Max}(\sO_K)}\frac{|V^*\otimes_{k[G]}{\rm Lie}(C^{\otimes n})(\sO_L/\mathfrak p\sO_L)|_{F[t]}}{|V^*\otimes_{k[G]}C^{\otimes n}(\sO_L/\mathfrak p\sO_L)|_{F[t]}}\in 1+t^{-1}F[[t^{-1}]].\end{eqnarray*}
\end{lem}
Since $G$ is arbitrary, then $V^*\otimes_{k[G]}\sO_L$ may not be a locally free $\sO_K$-module. So we can't use the method of \cite{F} to study $L(n,\,\rho)$. However, ${\rm Hom}_{k[G]}(V,\,\sO_L)$ is a locally free $\sO_K$-module. In order to study $L(n,\,\rho)$, we must consider the following convergent infinite product
\begin{eqnarray*}
L(C^{\otimes n},\rho)=\prod_{\mathfrak p\in{\rm Max}(\sO_K)}\frac{|{\rm Hom}_{k[G]}(V,\,{\rm Lie}(C^{\otimes n})(\sO_L/\mathfrak p\sO_L))|_{F[t]}}{|{\rm Hom}_{k[G]}(V,\,C^{\otimes n}(\sO_L/\mathfrak p\sO_L))|_{F[t]}}.
\end{eqnarray*}
In this paper, we also prove the following lemma.
\begin{lem}\label{3241}
We have
\begin{eqnarray*}
L(C^{\otimes n},\,\rho)=\prod_{\mathfrak p\in{\rm Max}(\sO_K)}{\rm det}_{F[[t^{-1}]]}\Big(1-\frac{\rho({\rm Frob}_\mathfrak P)}{|\kappa_\mathfrak p|^n_{k[t]}},\;V_{I_\mathfrak P}[[t^{-1}]]\Big)^{-1}\in 1+t^{-1}F[[t^{-1}]]
\end{eqnarray*}
and
$$\frac{L(n,\,\rho)}{L(C^{\otimes n},\,\rho)}\in F(t)^\times.$$
This means that $L(n,\,\rho)$ and $L(C^{\otimes n},\,\rho)$ are equal up to a rational polynomial.
\end{lem}

For any matrix $(a_{ij})$ over a $k$-algebra, denote by $(a_{ij})^{(q^s)}=(a_{ij}^{q^s})$. Let $M_n(\sO_K)$ be the ring of
$n\times n$-matrices over $\sO_K$ and let $M_n(\sO_K)\{\tau\}$ be the ring over $M_n(\sO_K)$ generated by $\tau$ with the relation $\tau P=P^{(q)}\tau$ for any $P\in M_n(\sO_K)$.
\begin{defn}\label{3134}
An abelian $t$-module $E$ over $\sO_K$ is a $k[t]$-module scheme $E$ over $\sO_K$ whose underlying $k$-vector space scheme is
isomorphic to $\textbf G_{\rm a}^n$ for some positive integer $n$ and the $k[t]$-module
structure on $E$ is given by
$$E:k[t]\to{\rm End}_{k\hbox{-}{\rm group}}(E)=M_n(\sO_K)\{\tau\}$$
such that $E(t)=\sum_{s=0}^rA_s\tau^s$
for some $A_0,\ldots,A_r\in M_n(\sO_K)$ with $(A_0-zI_n)^n=0$. The integer $n$ is called the dimension of $E$.
\end{defn}
\begin{rem}
An abelian $t$-module $E$ over $\sO_K$ defines a functor
$$E:\{\sO_K\hbox{-algebras}\}\to\{k[t]\hbox{-modules}\}.$$
The $n$-th tensor power $C^{\otimes n}$ of Carlitz module $C$ is an abelian $t$-module of dimension $n$.
\end{rem}

\begin{defn}
Define the $\infty$-adic $L$-value $L(E,\,\rho)$ of $E$ twisted by $\rho$ to be the convergent infinite product
$$L(E,\,\rho)=\prod_{\mathfrak p\in{\rm Max}(\sO_K)}\frac{|{\rm Hom}_{k[G]}(V,\,{\rm Lie}(E)(\sO_L/\mathfrak p\sO_L))|_{F[t]}}{|{\rm Hom}_{k[G]}(V,\,E(\sO_L/\mathfrak p\sO_L))|_{F[t]}}\in 1+t^{-1}F[[t^{-1}]].$$
\end{defn}

Let $K_\infty=K\otimes_{k(t)}k((t^{-1}))$ and $L_\infty=L\otimes_{k(t)}k((t^{-1}))$. There exists a unique power series
$$\exp_EX=\sum_{s\geq0}e_sX^{(q^s)}$$
with $X=(X_1,\ldots,X_n)^T$ and $e_s\in M_n(K_\infty)$ such that $e_0=I_n$ and
$$\exp_E(A_0X)=\sum_{s=0}^rA_s(\exp_EX)^{(q^s)}.$$
We get a continuous and open $k[t][G]$-linear map
$$\exp_E:{\rm Lie}(E)(L_\infty)\to E(L_\infty).$$
By \cite[Lemma 1.7]{F}, we get a $k((t^{-1}))[G]$-module structure on ${\rm Lie}(E)(L_\infty)$. By normal basis theorem, $L$ is a free $K[G]$-module of rank one. Then
${\rm Lie}(E)(L_\infty)$ is a free $k((t^{-1}))[G]$-module of rank $[K:k(t)]\dim E$.

\begin{defn}
Suppose $S$ is a commutative ring. A perfect complex of $S$-modules is a bounded complex of projective $S$-modules of finite type.
 Let $D^{\rm per}(S)$ be the full subcategory of $D(S)$ consisting of all objects which can be
 represented by a perfect complex. For any perfect complex $C=(C^i)_{i\in\Z}$, the determinant $\det_S(C)$ of $C$ is defined by $\bigotimes_{i\in\Z}\det_S(C^i)^{(-1)^i}$.
%For any two morphism $f,\,g:C\to D$, we have a canonical isomorphism
%$$\det(C^\bullet\x{f}D^\bullet)\simeq\det(C^\bullet)\otimes\det(D^\bullet)^{-1}\simeq\det(C^{\bullet}\x{g}\det(D^\bullet).$$
For any two isomorphisms $f,\,g:\det_S(C_1)\simeq \det_S(C_2)$ for some $C_1,\,C_2\in D^{\rm perf}(S)$, denote by $[f:g]$
the element of $S^\times$ such that $f=[f:g]g$.
\end{defn}

\begin{defn}\label{16} Let $H$ be a finite abelian group of order prime to $p$. Let $V$ be a free $F((t^{-1}))[H]$-module of finite rank.

(1) A lattice in $V$ is a finitely generated projective sub-$F[t][H]$-module $\Lambda$ of $V$ such that
the natural morphism $\Lambda\otimes_{F[t][H]}F((t^{-1}))[H]\to V$ is isomorphic.

(2) Let $\Lambda_1$ and $\Lambda_2$ be two lattices of $V$.
Any isomorphism $\det_{F[t][H]}(\Lambda_1)\simeq\det_{F[t][H]}(\Lambda_2)$ defines an isomorphism $\det_{F((t^{-1}))[H]}\Big(\Lambda_1\otimes_{F[t][H]}F((t^{-1}))[H]\Big)\simeq
\det_{F((t^{-1}))[H]}\Big(\Lambda_2\otimes_{F[t][H]}F((t^{-1}))[H]\Big).$ Let $f$ be the composite
\begin{eqnarray*}
&&{\rm det}_{F((t^{-1}))[H]}(V)\simeq{\rm det}_{F((t^{-1}))[H]}\Big(\Lambda_1\otimes_{F[t][H]}F((t^{-1}))[H]\Big)\\&\simeq&
{\rm det}_{F((t^{-1}))[H]}\Big(\Lambda_2\otimes_{F[t][H]}F((t^{-1}))[H]\Big)\simeq{\rm det}_{F((t^{-1}))[H]}(V).
\end{eqnarray*}
Then the image of $[f:1]$ in $F((t^{-1}))[H]^\times/F[H]^\times$ does not depend on the choice of the isomorphism $\det_{F[t][H]}(\Lambda_1)\simeq\det_{F[t][H]}(\Lambda_2)$.
Denote by $[\Lambda_1:\Lambda_2]_{F((t^{-1}))[H]}$ the unique monic representative
of $[f:1]$ in $F((t^{-1}))[H]^\times$.
\end{defn}
Define
$$H(E,\,\rho)={\rm coker}\Big({\rm Hom}_{k[G]}(V,\,{\rm Lie}(E)(L_\infty))\x{\exp_E}\frac{{\rm Hom}_{k[G]}(V,\,E(L_\infty))}{{\rm Hom}_{k[G]}(V,\,E(\sO_L))}\Big).$$
\begin{thm}\label{391}
Let $E$ be an abelian $t$-module over $\sO_K$. Let $\rho:G\to{\rm GL}_F(V)$ be a finite dimensional $F$-linear representation of the Galois group of a finite Galois extension $L$ of $K$. Then ${\rm Hom}_{k[G]}(V,\,{\rm Lie}(E)(\sO_L))$ and ${\rm Hom}_{k[G]}(V,\,\exp_E^{-1}(E(\sO_L)))$ are lattices of the finite dimensional $F((t^{-1}))$-vector space ${\rm Hom}_{k[G]}(V,\,{\rm Lie}(E)(L_\infty))$ and $H(E,\,\rho)$ is a finite $F[t]$-module.
Moreover,
\begin{eqnarray*}
L(E,\,\rho)=[{\rm Hom}_{k[G]}(V,\,{\rm Lie}(E)(\sO_L)):{\rm Hom}_{k[G]}(V,\,\exp_E^{-1}(E(\sO_L)))]_{F((t^{-1}))}\cdot|H(E,\,\rho)|_{F[t]}.
\end{eqnarray*}
\end{thm}
It is convenient to group all $L(E,\,\rho)$ together in an equivariant $L$-value $L(E,\,G)\in 1+t^{-1}k[[t^{-1}]][G]$ such that
\begin{eqnarray}\label{3165}
L(E,\,\rho)={\rm det}_{F((t^{-1}))}\Big(L(E,\,G),\,V((t^{-1}))\Big)\in 1+t^{-1}F[[t^{-1}]].
\end{eqnarray}
If $k[G]$ is not commutative or regular, we can't talk about the determinant for any $k[G]$-modules. So we must assume that $G$ is an abelian group of order prime to $p$ for the existence of an equivariant $L$-value $L(E,\,G)\in 1+ t^{-1}k[[t^{-1}]][G]$. So in the rest of this section, we assume that $G$ is an abelian group of order prime to $p$.

\begin{defn}
For an abelian $t$-module $E$ over $\sO_K$, the equivariant $\infty$-adic $L$-value $L(E,\,G)$ of $E$ on $L/K$ is defined to be the convergent infinite product
$$L(E,\,G)=\prod_{\mathfrak p\in{\rm Max}(\sO_K)}\frac{|{\rm Lie}(E)(\sO_L/\mathfrak p\sO_L)|_{k[t][G]}}{|E(\sO_L/\mathfrak p\sO_L)|_{k[t][G]}}\in k((t^{-1}))[G]^\times.$$
Define
$$H(E,\,G)=\frac{E(L_\infty)}{\exp_E({\rm Lie}(E)(L_\infty))+E(\sO_L)}.$$
\end{defn}
In this paper, we also prove the following theorem.
\begin{thm}\label{18}
Let $E$ be an abelian $t$-module over $\sO_K$. Let $L$ be a finite Galois extension of $K$ whose Galois group $G$ is an abelian group of order prime to $p$.

(1) Then ${\rm Lie}(E)(\sO_L)$ and $\exp_E^{-1}(E(\sO_L))$ are lattices of the free $k((t^{-1}))[G]$-module ${\rm Lie}(E)(L_\infty)$ of finite rank and $H(E,\,G)$
is a finite $k[t][G]$-module. We have
$$L(E,\;G)=[{\rm Lie}(E)(\sO_L):\exp_E^{-1}(E(\sO_L))]_{k((t^{-1}))[G]}\cdot|H(E,\;G)|_{k[t][G]}\in k((t^{-1}))[G]^\times.$$

(2) (\ref{3165}) holds for any finite dimensional $F$-linear representation $\rho:G\to{\rm GL}_F(V)$.

\end{thm}
\begin{rem}
These theorems generalize the class number formula of \cite{T} for Drinfeld modules and that of \cite{F} for abelian $t$-modules and the equivariant class number formula of \cite{AT} for cyclotomic extension of function fields.
\end{rem}

The paper is organized as follows. In section 2, we construct an equivariant version of Anderson's trace formula. In section 3, we construct the theory of equivariant crystals and prove the trace formula of equivariant $L$-functions. In section 4, we express the $v$-adic equivariant special values of shtukas in terms of the determinant of the extension groups of shtukas under some local analytic conditions. In section 5, we prove an equivariant class number formula for abelian $t$-modules. In section 6, we calculate the special values of Artin $L$-functions for Galois representations.

\textbf{Acknowledgements.}
I would like to thank David Goss and Lenny Taelman for their continuous encouragement.
My research is supported by the NSFC grant no. 00713210010018.

\section{Equivariant version of Anderson's trace formula}
In this section, let $G$ be a finite abelian group of order prime to $p$.
Let $R$ and $A$ be two $k$-algebras. For any $R$-module $M$, let $M[G]=R[G]\otimes_R M$.
\begin{lem}\label{zhong}
Let $M$ be an $R[G]$-module. Define
\begin{eqnarray*}
\pi:M[G]\to M,&\,&\bigoplus_{g\in G}g\otimes m_g\mapsto\sum_{g\in G}gm_g;\\
i:M\to M[G],&\,&m\mapsto\frac{1}{|G|}\sum_{g\in G}g\otimes g^{-1}m.
\end{eqnarray*}
Then $\pi$ and $i$ are $R[G]$-linear and $\pi\circ i={\rm id}_M$. Then $M[G]=i(M)\oplus N$ for $N={\rm ker}(\pi)$.

(1) Then $M$ is a projective $R[G]$-module if and only if $M$ is projective as an $R$-module.

(2) One has $M\otimes_kA\oplus N\otimes_kA\simeq (M\oplus N)\otimes_kA\simeq M[G]\otimes_kA=M\otimes_kA[G].$
For any $\phi\in{\rm End}_{A[G]}(M\otimes_kA)$, $\phi\oplus 0:M\otimes_kA\oplus N\otimes_kA\to M\otimes_kA\oplus N\otimes_kA$ induces an $A[G]$-linear map $\Phi:M\otimes_kA[G]\to M\otimes_kA[G]$. If $\dim_kM<\infty$, define
$${\rm det}_{A[G][t]}(1-t\phi,\;M\otimes_kA[t])={\rm det}_{A[G][t]}(1-t\Phi,\;M\otimes_kA[G][t])\in 1+A[G][t].$$
\end{lem}

%\begin{exam}
%The $n$-th tensor power of Carlitz module is an abelian $t$-module and $\zeta(R,\,n)=L(C^{\otimes n}/R)$.
%\end{exam}

%\begin{rem}
%Suppose the leading coefficient $A_r$ is invertible in $M_n(R)$, then $L(E/R)$ is the value at $s=0$ of the Goss L-function $L(E/R,\,s)$. The zeta value
%$\zeta(R,\,n)$ is the value at $s=n-1$ of the Goss L-function $L(C/R)$ for the Carlitz-module $C$.
%\end{rem}
\begin{defn}
Let $M$ be a $k[G]$-module and let $\phi\in{\rm End}_{A[G]}(M\otimes_kA)$.
A finitely generated $k[G]$-module $M_0$ of $M$ is called a nucleus of $\phi$ if there exists an exhaustive increasing filtration of $M$ by finitely generated $k[G]$-submodules $M_0\subset M_1\subset M_2\subset\cdots$ such that $\phi(M_{i+1}\otimes_kA)\subset M_i\otimes_kA$ for any $i\geq0$. One can show that ${\rm det}_{A[G][t]}(1-t\phi,M_0\otimes_kA[t])$ does not depend on the choice of the nucleus $M_0$ of $\phi$.
Define
$${\rm det}_{A[G][t]}(1-t\phi,M\otimes_kA[t])={\rm det}_{A[G][t]}(1-t\phi,\;M_0\otimes_kA[t]).$$

\end{defn}
\begin{defn}
Let $M$ and $N$ be two $R[G]$-modules. For any $n\geq1$, an $A[G]$-linear map $\phi:M\otimes_kA\to N\otimes_kA$ is called $n$-th Frobenius (resp. $n$-th Cartier) over $A$ on $R$ if
$$\phi(rm)=r^{q^n}\phi(m) \;\;(\hbox{resp. }\phi(r^{q^n}m)=r\phi(m))\hbox{ for any }r\in R\hbox{ and }m\in M\otimes_kA.$$
%In particular, if $M=N$, we call $\phi$ is a $n$-th Frobenius (resp. $n$-th Cartier) operator on $M$.
\end{defn}

\begin{defn}\label{yuan}
Let $R$ be a finitely generated $k$-algebra and $M$ a finitely generated $R[G]$-module. For any $n\geq1$, let $C_n:M\otimes_kA\to M\otimes_kA$ be an $n$-th Cartier operator over $A$ on $R$. By \cite[Proposition 6]{A1}, $C_1,C_2,\ldots,C_n$ has a common nucleus $M_n$ for any $n\geq1$. Then
$${\rm det}_{A[G][t]}\Big(1-\sum_{i=1}^nt^iC_i,\;M_n\otimes_kA[t]\Big)\in 1+tA[G][t]$$ does not depend on the choice of $M_n$.
Define
$${\rm det}_{A[G][[t]]}\Big(1-\sum_{n=1}^\infty t^nC_n,\;M\otimes_kA[[t]]\Big)=\lim_{n\to\infty}{\rm det}_{A[G][t]}\Big(1-\sum_{i=1}^nt^iC_i,\;M_n\otimes_kA[t]\Big)\in 1+tA[G][[t]].$$
\end{defn}

\begin{thm}\label{ji}
Let $R$ be a finitely generated regular domain over $k$ of Krull dimension $r$. Let $\Omega_R=\wedge^r\Omega_{R/k}^1$. By \cite[2.6]{A1}, we have a Cartier operator $C$ on $\Omega_R$. Let $M$ be a finitely generated projective $R[G]$-module. For each $n\geq1$, let $\tau_n:M\otimes_kA\to M\otimes_kA$ be an $n$-th Frobenius operator over $A$ on $R$. Let $\widehat M={\rm Hom}_R(M,\,\Omega_R)$. Define the adjoint operator $C_n$ on $\widehat M\otimes_kA\simeq{\rm Hom}_{R\otimes_kA}(M\otimes_kA,\,\Omega_R\otimes_kA)$ of $\tau_n$ by the rule
$$(C_n\hat m)(m)=(C^{n[k:\F_p]}\otimes{\rm id}_A)(\hat m(\tau_n(m)))\hbox{ for any }m\in M\otimes_kA\hbox{ and }\hat m\in\widehat M\otimes_kA.$$
Let $I$ be a maximal ideal of $R$ of degree $d=[R/I:k]$. Then
$${\rm det}_{A[G][[t]]}\Big(1-\sum_{n=1}^\infty t^n\tau_n,\;M/IM\otimes_kA[[t]]\Big)^{-1}\in 1+t^dA[G][[t^d]].$$
We have
$$\prod_{I\in{\rm Max}(R)}{\rm det}_{A[G][[t]]}\Big(1-\sum_{n=1}^\infty t^n\tau_n,\;M/IM\otimes_kA[[t]]\Big)^{-1}={\rm det}_{A[G][[t]]}\Big(1-\sum_{n=1}^\infty t^nC_n,\,\widehat M\otimes_kA[[t]]\Big)^{(-1)^{r-1}}.$$
\end{thm}
\begin{proof}
By Lemma \ref{zhong}, we may assume $G=\{1\}.$ For the sequence $\{\tau_n\}_{n\geq1}$, there exists a unique $n$-th Frobenius operator $\tau_n'$ on $M\otimes_kA$ over $A$ on $R$ for each $n$ such that
$$1-\sum_{n=1}^\infty t^n\tau_n=(1-t\tau_1')(1-t^2\tau_2')(1-t^3\tau_3')\cdots.$$
Let $C_n':\widehat M\otimes_kA\to\widehat M\otimes_kA$ be the adjoint of $\tau'_n$. Then we have
$$1-\sum_{n=1}^\infty t^nC_n=\cdots(1-t^3C_3')(1-t^2C_2')(1-tC_1').$$

For any $n\geq1$, let $k_n$ be the extension field of $k$ of degree $n$. For any object $\sF$ over $k$, let $\sF_n=k_n\otimes_k\sF$. For any $k$-linear map $\phi:\sF\to\sF$, denote the $k_n$-linear map ${\rm id}_{k_n}\otimes\phi:\sF_n\to\sF_n$ also by $\phi$. Then $\tau_n':M_n\otimes_{k_n}A_n\to M_n\otimes_{k_n}A_n$ is a 1-th Frobenius operator over $A_n$ on $R_n$. By \cite[Theorem 1]{A1}, we have
$$\prod_{J\in{\rm Max}(R_n)}{\rm det}_{A_n[[t]]}\Big(1- t\tau_n',\;M_n/JM_n\otimes_{k_n}A_n[[t]]\Big)^{-1}={\rm det}_{A_n[[t]]}\Big(1- tC_n',\,\widehat M_n\otimes_{k_n}A_n[[t]]\Big)^{(-1)^{r-1}}.$$
For any maximal ideal $I$ of $R$ of degree $d$, $k_n \otimes_k R/I\simeq\prod_{I\subset J\in{\rm Max}(R_n)}R_n/J$ and $[R_n/J:k_n]=\frac{d}{(n,\,d)}$. By \cite[Lemma 8.1.4]{BP}, we have
\begin{eqnarray*}
&&{\rm det}_{A[[t]]}\Big(1-t\tau_n',\,M/IM\otimes_kA[[t]]\Big)\\&=&{\rm det}_{A_n[[t]]}\Big(1-t\tau_n',\,M_n/IM_n\otimes_{k_n}A_n[[t]]\Big)\\
&=&\prod_{I\subset J\in{\rm Max}(R_n)}{\rm det}_{A_n[[t]]}\Big(1-t\tau_n',\,M_n/JM_n\otimes_{k_n}A_n[[t]]\Big)\in 1+t^{\frac{d}{(n,\,d)}}A[[t^{\frac{d}{(n,\,d)}}]].
\end{eqnarray*}
Substituting $t$ by $t^n$,
we get
\begin{eqnarray*}
&&{\rm det}_{A[[t]]}\Big(1-\sum_{n=1}^\infty t^n\tau_n,\;M/IM\otimes_kA[[t]]\Big)^{-1}\\
&=&\prod_{n=1}^\infty {\rm det}_{A[[t]]}\Big(1- t^n\tau_n',\;M/IM\otimes_kA[[t]]\Big)^{-1}\in 1+t^{\frac{nd}{(n,\,d)}}A[[t^{\frac{nd}{(n,\,d)}}]]\subset 1+t^dA[[t^d]]
\end{eqnarray*} and
\begin{eqnarray*}
&&\prod_{I\in{\rm Max}(R)}{\rm det}_{A[[t]]}\Big(1-t^n\tau_n',\;M/IM\otimes_kA[[t]]\Big)^{-1}\\
&=&\prod_{I\in{\rm Max}(R)}\prod_{I\subset J\in{\rm Max}(R_n)}{\rm det}_{A_n[[t]]}\Big(1-t^n\tau_n',\;M_n/JM_n\otimes_{k_n}A_n[[t]]\Big)^{-1}\\
&=&\prod_{J\in{\rm Max}(R_n)}{\rm det}_{A_n[[t]]}\Big(1-t^n\tau_n',\;M_n/JM_n\otimes_{k_n}A_n[[t]]\Big)^{-1}\\
&=&{\rm det}_{A_n[[t]]}\Big(1- t^nC_n',\,\widehat M_n\otimes_{k_n}A_n[[t]]\Big)^{(-1)^{r-1}}\\
&=&{\rm det}_{A[[t]]}\Big(1- t^nC_n',\,\widehat M\otimes_{k}A[[t]]\Big)^{(-1)^{r-1}}.
\end{eqnarray*}

Then we have
\begin{eqnarray*}
&&\prod_{I\in{\rm Max}(R)}{\rm det}_{A[[t]]}\Big(1-\sum_{n=1}^\infty t^n\tau_n,\;M/IM\otimes_kA[[t]]\Big)^{-1}\\
&=&\prod_{I\in{\rm Max}(R)}\prod_{n=1}^\infty{\rm det}_{A[[t]]}\Big(1- t^n\tau_n',\;M/IM\otimes_kA[[t]]\Big)^{-1}\\
&=&\prod_{n=1}^\infty\prod_{I\in{\rm Max}(R)}{\rm det}_{A[[t]]}\Big(1- t^n\tau_n',\;M/IM\otimes_kA[[t]]\Big)^{-1}\\
&=&\prod_{n=1}^\infty{\rm det}_{A[[t]]}\Big(1-t^nC_n',\,\widehat M\otimes_kA[[t]]\Big)^{(-1)^{r-1}}\\
&=&{\rm det}_{A[[t]]}\Big(1-\sum_{n=1}^\infty t^nC_n,\,\widehat M\otimes_kA[[t]]\Big)^{(-1)^{r-1}}.
\end{eqnarray*}
\end{proof}

\section{Trace formula for equivariant $L$-functions of crystals}
In this section, suppose $G$ is a finite group and $A$ is a $k$-algebra. All tensors are over $k$ and all schemes are finite type over $k$.
For any scheme $X$, let $\sigma_X:X\to X$ be the morphism induced by the $q$-th power map on $\sO_X$.
\begin{defn}\label{35}
(1) A coherent $\tau$-sheaf over $A$ on $X$ is a pair $(\sF,\,\tau_\sF)$ consisting of a coherent $\sO_X$-module $\sF$ and an $\sO_X\otimes A$-linear homomorphism
$$\tau_\sF:\bigoplus_{1\leq n\ll\infty}(\sigma_X^n)^*\sF\otimes A\to\sF\otimes A.$$

(2) A homomorphism $(\sF,\,\tau_\sF)\to(\sG,\,\tau_\sG)$ is an $\sO_X$-linear map $\phi:\sF\to\sG$ such that the following diagram commutes
\[\xymatrix{\bigoplus_{1\leq n\ll\infty}(\sigma_X^n)^*\sF\otimes A\ar[rr]^{\tau_\sF}\ar[d]^{\oplus(\sigma_X^n)^*\phi\otimes {\rm id}_A}&&\sF\otimes A\ar[d]^{\phi\otimes {\rm id}_A}\\
\bigoplus_{1\leq n\ll\infty}(\sigma_X^n)^*\sG\otimes A\ar[rr]^{\tau_\sG}&&\sG\otimes A.}\]
\end{defn}

(3) A coherent $\tau$-sheaf $(\sF,\,\tau_\sF)$ is called nilpotent if there exists a decreasing filtration $\{\sF_m\}_{m\in\N}$ of $\sF$ by coherent submodules such that $\sF_0=\sF$, $\sF_m=0$ for $m\gg0$ and $$\tau_\sF\Big(\bigoplus_{1\leq n\ll\infty}(\sigma_X^n)^*\sF_m\otimes A\Big)\subset\sF_{m+1}\otimes A\hbox{ for any }m.$$

\begin{defn}
Let $\sA$ be an abelian category.

(1) A $G$-module in $\sA$ is a pair $(\sF,\,\alpha)$ consisting of an object $\sF$ in $\sA$ and a homomorphism $\alpha:G\to{\rm Aut}(\sF)$ of groups.

(2) Let $K_0(\sA)$ be the Grothendieck group of $\sA$ generated by the isomorphic classes $[M]$ of objects in $\sA$, modulo relations $[M_2]=[M_1]+[M_3]$ for every short exact sequence in $\sA$
$$0\to M_1\to M_2\to M_3\to0.$$
\end{defn}
\begin{defn}\label{3201}
(1) Let ${\rm Coh}_\tau(X,\,A)$ be the category of coherent $\tau$-sheaves over $A$ on $X$ and let ${\rm NilCoh}_\tau(X,\,A)$ be its full subcategory consisting of nilpotent objects.

(2) Let ${\rm Coh}^G_\tau(X,\,A)$ and ${\rm NilCoh}^G_\tau(X,\,A)$ be the category of $G$-modules in ${\rm Coh}_\tau(X,\,A)$ and ${\rm NilCoh}_\tau(X,\,A)$, respectively.
%Let ${\rm Crys}^G(X,\,A)=\frac{{\rm Coh}^G_\tau(X,\,A)}{{\rm NilCoh}^G_\tau(X,\,A)}.$

(3) A morphism $\alpha$ in ${\rm Coh}^G_\tau(X,\,A)$ is called nil-isomorphism if ${\rm ker}(\alpha)$ and ${\rm coker}(\alpha)$ are nilpotent.
\end{defn}

\begin{prop}
The category ${\rm Coh}^G_\tau(X,\,A)$ is abelian and ${\rm NilCoh}^G_\tau(X,\,A)$ is its abelian subcategory stable under subquotient and extension. By \cite [Proposition 2.3.4]{BP}, the quotient category ${\rm Crys}^G(X,\,A)$ of ${\rm Coh}^G_\tau(X,\,A)$ by ${\rm NilCoh}^G_\tau(X,\,A)$ is also abelian.
Let $q_X:{\rm Coh}^G_\tau(X,\,A)\to{\rm Crys}^G(X,\,A)$ be the localization map.  Then $\alpha$ is a nil-isomorphism in ${\rm Coh}^G_\tau(X,\,A)$ if and only if $q_X(\alpha)$ is an isomorphism in ${\rm Crys}^G(X,\,A)$. Moreover, $q_X$ is exact and induces a homomorphism
$$q_X:K_0({\rm Coh}^G_\tau(X,\,A))\to K_0({\rm Crys}^G(X,\,A)).$$
\end{prop}

\begin{rem}
Definition \ref{35} and Definition \ref{3201} of $\tau$-sheaves and crystals are stronger than that of \cite{BP}. Under the definition of crystals in \cite{BP}, one only get a cohomological trace formula of $L$-functions up to unipotent polynomials. However, we get an explicit cohomological trace formula of equivariant $L$-functions of crystals in our setting.
\end{rem}

\subsection{Inverse image}
Let $f:X\to Y$ be a $k$-morphism of schemes and let $i\geq0$. For any coherent $\tau$-sheaf $\underline\sG=(\sG,\tau_\sG)$ over $A$ on $Y$, denote by $L_if^*\underline\sG=(L_if^*\sG,\,\tau_{L_if^*\sG})$ where $\tau_{L_if^*\sG}$ is the composition
\begin{eqnarray*}
\bigoplus_{1\leq n\ll\infty}(\sigma_X^n)^*L_if^*\sG\otimes A\x{{\rm can.}}\bigoplus_{1\leq n\ll\infty} L_if^*(\sigma_Y^n)^*\sG\otimes A\x{L_if^*\tau_\sG}L_if^*\sG\otimes A.
\end{eqnarray*}
If $\underline\sG$ is a $G$-module in ${\rm Coh}_\tau(Y,\,A)$, then we get a $G$-module $L_if^*\underline\sG$ in ${\rm Coh}_\tau(X,\,A)$.
Thus we get a functor
$$L_if^*:{\rm Coh}_\tau^G(Y,\,A)\to{\rm Coh}_\tau^G(X,\,A).$$

\begin{lem}\label{zhao}
Let $\underline\sG\in{\rm Coh}_\tau(X,\,A)$. For any $i\geq 1$, $L_if^*\underline\sG$ is a nilpotent coherent $\tau$-sheaf on $X$.
\end{lem}
\begin{proof}
Factors $f=hg$ with a closed immersion $g$ and a flat morphism $h$. Then
$L_if^*\sG=L_ig^*(h^*\sG)$ and $L_if^*\underline\sG=L_ig^*(h^*\underline\sG)$. Hence we may assume that $f$ is a closed immersion defined by an ideal $\sI$ of $\sO_Y$.
Define a decreasing filtration $\{\sF_m\}_{m\geq0}$ of $L_if^*\sG={\rm Tor}_i^{\sO_Y}(\sG,\,\sO_Y/\sI)$ by
$$\sF_m={\rm im}\Big({\rm Tor}_i^{\sO_Y}(\sG,\,\sO_Y/\sI^{m+1})\to{\rm Tor}_i^{\sO_Y}(\sG,\,\sO_Y/\sI)\Big).$$
It remains to show that $\sF_m=0$ for $m\gg0$ and $$\tau_{L_if^*\sG}\Big(\bigoplus_{1\leq n\ll\infty}(\sigma_X^n)^*\sF_m\otimes A\Big)\subset \sF_{m+1}\otimes A.$$
The problem is local on $Y$. So we may assume $Y$ is affine. Then we can choose an exact sequence
$$0\to(\sK,\tau_\sK)\to(\sP_i,\tau_{\sP_i})\to\cdots\to(\sP_1,\tau_{\sP_1})\to(\sG,\,\tau_\sG)\to0 $$
of coherent $\tau$-sheaves over $A$ on $Y$ such that $\sP_1,\ldots, \sP_i$ are free $\sO_Y$-modules. Then
$${\rm Tor}_i^{\sO_Y}(\sG,\,\sO_Y/\sI^{m+1})=\sK\cap \sI^{m+1}\sP_i/\sI^{m+1}\sK$$
and hence
$$ \sF_m=(\sK\cap \sI^{m+1}\sP_i)+\sI\sK/\sI\sK.$$
By Artin-Rees lemma, there exists $m_0\geq0$ such that for any $m\geq m_0$, we have
$$\sK\cap \sI^{m+1}\sP_i=\sI^{m+1-m_0}(\sK\cap \sI^{m_0}\sP_i). $$
We have
$\sK\cap \sI^{m+1}\sP_i\subset \sI\sK $ and hence $\sF_m=0$ for any $m\geq m_0$.
This shows
$$\tau_{\sP_i}\Big(\bigoplus_{1\leq n\ll\infty} \sK\cap \sI^{m+1}\sP_i\otimes A\Big)\subset \sK\cap \sI^{q(m+1)}\sP_i\otimes A$$
for any $m\geq0$. Then for any $m\geq0$, we have $$\tau_{L_if^*\sG}\Big(\bigoplus_{1\leq n\ll\infty} (\sigma^n_{X})^*\sF_m\otimes A\Big)\subset \sF_{m+1}\otimes A.$$
\end{proof}

\begin{lem}\label{nong}
There is a unique functor $f^*:{\rm Crys}^G(Y,\,A)\to{\rm Crys}^G(X,\,A)$ such that the diagram
\[\xymatrix{{\rm Coh}^G_\tau(Y,\,A)\ar[r]^{f^*}\ar[d]^{q_Y}&{\rm Coh}^G_\tau(X,\,A)\ar[d]^{q_X}\\{\rm Crys}^G(Y,\,A)\ar[r]^{f^*}&{\rm Crys}^G(X,\,A)}\]
commutes. Moreover, the functor $f^*:{\rm Crys}^G(Y,\,A)\to{\rm Crys}^G(X,\,A)$ is exact and hence it induces a homomorphism
$$f^*:K_0({\rm Crys}^G(Y,\,A))\to K_0({\rm Crys}^G(X,\,A)).$$
\end{lem}
\begin{proof}
If $\alpha$ is a nil-isomorphism in ${\rm Coh}^G_\tau(Y,\,A)$, then $f^*(\alpha)$ is a nil-isomorphism in ${\rm Coh}^G_\tau(X,\,A)$ and $q_Xf^*(\alpha)$ is an isomorphism.
Then the existence of $f^*:{\rm Crys}^G(Y,\,A)\to{\rm Crys}^G(X,\,A)$ holds by \cite[Theorem 2.2.3 (b)]{BP}.
Any short exact sequence in ${\rm Crys}^G(Y,\,A)$ is isomorphic to the image of a short exact sequence $0\to\underline\sF\to\underline\sG\to\underline\sH\to0$ in ${\rm Coh}^G_\tau(Y,\,A)$.
We have an exact sequence
$$L_1f^*\underline\sH\to f^*\underline\sF\to f^*\underline\sG\to f^*\underline\sH\to0$$
 in ${\rm Coh}^G_\tau(X,\,A)$. By Lemma \ref{zhao}, $L_1f^*\underline\sH\in{\rm NilCoh}_\tau^G(X,\,A)$. Then
$$0\to f^*\underline\sF\to f^*\underline\sG\to f^*\underline\sH\to0$$ is exact in ${\rm Crys}^G(X,\,A)$.
\end{proof}
\subsection{Galois cohomology of $G$-modules of crystals}
Let $\underline\sF=(\sF,\tau_\sF)\in{\rm Coh}_\tau^G(X,\,A)$. For any $i\geq0$, the $i$-th Galois cohomology $H^i(G,\,\sF)$ of $\sF$ is a coherent sheaf on $X$ and $\tau_\sF$ induces a map
$$\tau_{H^i(G,\,\sF)}:\bigoplus_{1\leq n\ll\infty}(\sigma_X^n)^*H^i(G,\,\sF)\otimes A\to H^i(G,\,\sF)\otimes A.$$
Thus we get a coherent $\tau$-sheaf $H^i(G,\,\underline\sF)=(H^i(G,\,\sF),\,\tau_{H^i(G,\,\sF)})$ over $A$ on $X$.
So we get two functors
\begin{eqnarray*}
&&H^i(G,\,\bullet):{\rm Coh}_\tau^G(X,\,A)\to{\rm Coh}_\tau(X,\,A);\\
&&H^i(G,\,\bullet):{\rm Crys}^G(X,\,A)\to{\rm Crys}(X,\,A).
\end{eqnarray*}
\begin{thm}\label{3155}
For any morphism $f:X\to Y$ of $k$-schemes, we have a commutative diagram
\[\xymatrix{{\rm Crys}^G(Y,\,A)\ar[r]^{f^*}\ar[d]^{H^i(G,\,\bullet)}&{\rm Crys}^G(X,\,A)\ar[d]^{H^i(G,\,\bullet)}\\{\rm Crys}(Y,\,A)\ar[r]^{f^*}&{\rm Crys}(X,\,A).}\]
\end{thm}
\begin{proof}
For any  $\underline\sF\in{\rm Coh}_\tau^G(Y,\,A)$, let $C^i(G,\,\sF)$ be the sheaf of sets of all maps from $G^i$ to $\sF$. Then $C^i(G,\,\sF)$ is a coherent sheaf on $X$. Define
$$d_i:C^i(G,\,\sF)\to C^{i+1}(G,\,\sF)$$
by
\begin{eqnarray*}
&&(d_if)(g_1,\ldots,g_{i+1})\\
&=&g_1f(g_2,\ldots,g_{i+1})+\sum_{j=1}^if(g_1,\ldots,g_{j-1},g_jg_{j+1},g_{j+2},\ldots,g_{i+1})+(-1)^{i+1}f(g_1,\ldots,g_i)
\end{eqnarray*}
for any $f\in C^i(G,\,\sF)$ and $(g_1,\ldots,g_{i+1})\in G^{i+1}$. By the definition of Galois cohomology, $H^i(G,\,\sF)$ is the $i$-th homology group of the complex
$$0\to C^0(G,\,\sF)\x{d_0}C^1(G,\,\sF)\to\cdots\to C^i(G,\,\sF)\x{d_i}C^{i+1}(G,\,\sF)\to\cdots.$$
Then
$$H^i(G,\,\sF)=H^iC^\bullet(G,\sF)\hbox{  and  }H^i(G,\,\underline\sF)=H^iC^\bullet(G,\,\underline\sF).$$
By Lemma \ref{nong}, we have
$$f^*H^i(G,\,\underline\sF)=f^*H^iC^\bullet(G,\,\underline\sF)\simeq H^i f^*C^\bullet(G,\,\underline\sF)=H^iC^\bullet(G,\,f^*\underline\sF)=H^i(G,\,f^*\underline\sF)\in{\rm Crys}(X,\,A).$$
\end{proof}
\subsection{Direct image for proper morphisms}
Let $f:X\to Y$ be a proper morphism and let $i\geq0$.
For any $\underline\sF\in{\rm Coh}_\tau(X,\,A)$, denote by $R^if_*\underline\sF$ the coherent sheaf $R^if_*\sF$ on $X$ together with the homomorphism
$$\tau_{R^if_*\sF}:\bigoplus_{1\leq n\ll\infty}(\sigma_Y^n)^*R^if_*\sF\otimes A\to R^if_*\sF\otimes A$$ induced by $\tau_\sF$. If $\underline\sF$ is a $G$-module in ${\rm Coh}_\tau(X,\,A)$, then $R^if_*\underline\sF$ is a $G$-module in ${\rm Coh}_\tau(Y,\,A)$.

\begin{thm}\label{351}
Let $\underline\sF=(\sF,\,\tau_\sF)\in {\rm Coh}_\tau(X,\,A)$. For any cartesian diagram
\begin{eqnarray}\label{ming}
\xymatrix{X'\ar[r]^{g'}\ar[d]^{f'}&X\ar[d]^f\\Y'\ar[r]^g&Y,}
\end{eqnarray}
of $k$-schemes, the natural homomorphism
$$g^*R^if_*\underline\sF\to R^if'_*g'^*\underline\sF$$
is a nil-isomorphism in ${\rm Coh}_\tau(Y',\,A)$.
\end{thm}
\begin{proof}
Factors $g=h\iota$ with a closed immersion $\iota$ and a flat morphism $h$. Consider the cartesian diagram
\[\xymatrix{X'\ar[r]^{\iota'}\ar[d]^{f'}&X''\ar[d]^{f''}\ar[r]^{h'}&X\ar[d]^f\\Y'\ar[r]^\iota&Y''\ar[r]^h&Y.}\]
By flat base change theorem, we have a natural isomorphism of $\sO_{Y''}$-modules
$$h^*R^if_*\sF\simeq R^if''_*h'^*\sF.$$
So we get a natural isomorphism $$h^*R^if_*\underline\sF\simeq R^if''_*h'^*\underline\sF$$
in ${\rm Coh}_\tau(Y'',\,A)$. The natural homomorphism $g^*R^if_*\underline\sF\to R^if'_*g'^*\underline\sF$ is the composite
$$g^*R^if_*\underline\sF\simeq \iota^*h^*R^if_*\underline\sF\simeq\iota^*R^if''_*(h'^*\underline\sF)\to R^if'_*\iota'^*(h'^*\underline\sF)\simeq R^if'_*g'^*\underline\sF.$$
Then we may assume that $g$ is a closed immersion defined by an ideal $\sI$ of $\sO_Y$.
For any $n\geq0$, let $\sF_n=\sF/I^{n+1}\sF.$ It remains to show the natural homomorphism
$$R^if_*\underline\sF\otimes_{\sO_Y}\sO_Y/\sI\to R^if'_*\underline{\sF_0}=R^if_*\underline{\sF_0}$$
is a nil-isomorphism in ${\rm Coh}_\tau(Y',\,A)$.
By \cite [Theorem 4.1.5 ] {G}, we have an isomorphism
$$\varprojlim_n(R^if_*\sF\otimes_{\sO_Y}\sO_Y/\sI^{n+1})\simeq\varprojlim_nR^if_*\sF_n.$$
Then we have a natural isomorphism
$$R^if_*\sF\otimes_{\sO_Y}\sO_Y/\sI\simeq (\varprojlim_nR^if_*\sF_n)\otimes_{\sO_Y}\sO_Y/\sI.$$
It remains to show that the natural map
$$\alpha:(\varprojlim_nR^if_*\underline\sF_n)\otimes_{\sO_Y}\sO_Y/\sI\to R^if_*\underline{\sF_0}$$
is a nil-isomorphism. By \cite [Theorem 4.1.5 ] {G}, $\{R^if_*\sF_n\}_{n\geq0}$ satisfies the Mittag-Leffler condition, thus
$${\rm im}\Big(\varprojlim_nR^if_*\sF_n\to R^if_*\sF_0\Big)={\rm im}\Big(R^if_*\sF_m\to R^if_*\sF_0\Big)\hbox{ for }m\gg0.$$
Since $\underline{\sF_m}\to\underline{\sF_0}$ is a nil-isomorphism, so is $R^if_*\underline{\sF_m}\to R^if_*\underline{\sF_0}$.
Then for $m\gg0$,
$${\rm cok}(\alpha)={\rm cok}\Big(\varprojlim_nR^if_*\underline\sF_n\to R^if_*\underline{\sF_0}\Big)={\rm cok}\Big(R^if_*\underline{\sF_m}\to R^if_*\underline{\sF_0}\Big)$$ is nilpotent.
Let
$$\sK_m=\ker\Big(R^if_*\sF\to R^if_*\sF_m\Big).$$
Also by \cite [Theorem 4.1.5 ] {G}, $\{\sK_m\}_{m\geq0}$ is an $\sI$-stable filtration of $R^if_*\sF$. There exists $m_0\in\N$ such that $\sK_m=\sI^{m-m_0}\sK_{m_0}$ for any $m\geq m_0$.
Then $\sK_m=\sI^{m-m_0}\sK_{m_0}\subset \sI\sK_0$ for $m>m_0$. Define a filtration $\{\sG_m\}_{m\geq0}$ of $\sK_0/\sI\sK_0$ by $\sG_m=\sI\sK_0+\sK_m/\sI\sK_0.$ Then $\sG_m=0$ for $m>m_0$. The short exact sequence
$$0\to\sI^{m+1}\sF\to\sF\to\sF_m\to0$$
shows that
$$\sK_m=\ker\Big(R^if_*\sF\to R^if_*\sF_m\Big)={\rm im}\Big(R^if_*\sI^{m+1}\sF\to R^if_*\sF_m\Big).$$
Then for any $m\geq0$, we have
$$\tau_{R^if_*\sF}\Big(\bigoplus_{1\leq n\ll\infty}(\sigma_Y^n)^*\sK_m\otimes A\Big)\subset\sK_{m+1}\otimes A.$$
Thus $\underline\sK_0/\sI\underline\sK_0\in{\rm NilCoh}_\tau(Y',\,A)$. The natural map
$${\sK_0}=\ker\Big(R^if_*\sF\to R^if_*{\sF_0}\Big)\to\ker\Big(R^if_*\sF\otimes_{\sO_Y}\sO_Y/\sI\x{\alpha}R^if_*{\sF_0}\Big)$$
is surjective, so is $\underline{\sK_0}/\sI\underline{\sK_0}\to\ker(\alpha)$. This proves $\ker(\alpha)$ is nilpotent, which completes the proof.
\end{proof}

\begin{cor}\label{lan}
Keep assumptions of Theorem \ref{351}. We have a homomorphism
$$Rf_*:K_0({\rm Crys}^G(X,\,A))\to K_0({\rm Crys}^G(Y,\,A)),\;[\underline\sF]\mapsto\sum_i(-1)^i[R^if_*\underline\sF]$$
of abelian groups and a commutative diagram
\[\xymatrix{K_0({\rm Crys}^G(X,\,A))\ar[r]^{Rf_*}\ar[d]^{g'^*}&K_0({\rm Crys}^G(Y,\,A))\ar[d]^{g^*}\\K_0({\rm Crys}^G(X',\,A))\ar[r]^{Rf'_*}&K_0({\rm Crys}^G(Y',\,A)).}\]
\end{cor}

\subsection{Extension by zero and proper base change}
\begin{lem}
Let $j:U\to X$ be an open immersion. Then the functor $j^*:{\rm Crys}^G(X,\,A)\to{\rm Crys}^G(U,\,A)$ has a left adjoint functor $j_!:{\rm Crys}^G(U,\,A)\to{\rm Crys}^G(X,\,A)$.
\end{lem}
\begin{proof}
By \cite[ Proposition 4.5.7 ]{BP}, $j^*: {\rm Crys}(X,\,A)\to{\rm Crys}(U,\,A)$ has a left adjoint functor $j_!:{\rm Crys}(U,\,A)\to{\rm Crys}(X,\,A)$. By \cite [Proposition 5.2 ]{T1}, ${\rm Crys}^G(U,\,A)$ and ${\rm Crys}^G(X,\,A)$ are equivalent to the categories of $G$-modules in ${\rm Crys}(U,\,A)$ and ${\rm Crys}(X,\,A)$, respectively. Then $j_!:{\rm Crys}(U,\,A)\to{\rm Crys}(X,\,A)$ defines a functor $j_!:{\rm Crys}^G(U,\,A)\to{\rm Crys}^G(X,\,A)$ which is left adjoint to $j^*:{\rm Crys}^G(X,\,A)\to{\rm Crys}^G(U,\,A)$.
\end{proof}

\begin{thm}\label{yun}
Let $f:X\to Y$ be a morphism of $k$-schemes. Then we can find a proper morphism $\bar f:\overline X\to Y$ and an open immersion $j:X\to\overline X$ such that $f=\bar f\circ j$. Define $Rf_!$ to be the composite
$$K_0({\rm Crys}^G(X,\,A))\x{j_!}K_0({\rm Crys}^G(\overline X,\,A))\x{R\bar f_*}K_0({\rm Crys}^G(Y,\,A)).$$
(1) The functor $Rf_!$ does not depend on the choice of $\bar f$ and $j$. Consider diagram (\ref{ming}).
We have a commutative diagram
\[\xymatrix{K_0({\rm Crys}^G(X,\,A))\ar[r]^{Rf_!}\ar[d]^{g'^*}&K_0({\rm Crys}^G(Y,\,A))\ar[d]^{g^*}\\K_0({\rm Crys}^G(X',\,A))\ar[r]^{Rf'_!}&K_0({\rm Crys}^G(Y',\,A)).}\]
(2) For any morphism $g:Y\to Z$ of $k$-schemes, we have a commutative diagram
\[\xymatrix{K_0({\rm Crys}^G(X,\,A))\ar[r]^{Rf_!}\ar[rd]^{R(gf)_!}&K_0({\rm Crys}^G(Y,\,A))\ar[d]^{Rg_!}\\&K_0({\rm Crys}^G(Z,\,A)).}\]
\end{thm}

\subsection{Trace formula for $L$-functions} In this subsection, $G$ is assumed to be a finite abelian group of order prime to $p$.

\begin{defn}Let $x={\rm Spec}\;k'$ for a finite field extension $k'/k$. Any $G$-module $\underline\sF=(\sF,\,\tau_\sF)$ in ${\rm Coh}_\tau(x,\,A)$ consisting of a finitely generated $k'[G]$-module $\sF$ together with an $n$-th Forbenius operator $\tau_n:\sF\otimes_kA\to\sF\otimes_kA$ over $A$ on $k'$ for finitely many positive integers $n$.
\label{3171}
The $L$-function of $\underline\sF\in{\rm Coh}_\tau^G(x,\,A)$ is defined to be
$$L^G(x,\,\underline\sF,\,t)={\rm det}_{A[G][t]}\Big(1-\sum_{1\leq n\ll\infty}t^n\tau_n,\,\sF\otimes_kA[t]\Big)^{-1}\in1+tA[G][[t]].$$

\end{defn}

\begin{defn}
Let $\underline\sF$ be a $G$-module in ${\rm Coh}_\tau(X,\,A)$ for a $k$-scheme $X$. For any closed point $x$ of $X$, let $i_x:{\rm Spec}\,\kappa(x)\to X$ be the closed immersion defined by $x$. By Theorem \ref{ji}, $L^G(x,\,i_x^*\underline\sF,\,t)\in 1+t^{d(x)}A[G][[t^{d(x)}]]$ where $d(x)=[\kappa(x):k]$. Then the product
$\prod_{x\in|X|}L^G(x,\,i_x^*\underline \sF,\,t)$ converges in $1+tA[G][[t]]$. Define the $L$-function of $\underline\sF\in{\rm Coh}_\tau(X,\,A)$ by
$$L^G(X,\,\underline\sF,\,t)=\prod_{x\in |X|}L^G(x,\,i_x^*\underline\sF,\,t)\in 1+tA[G][[t]].$$
\end{defn}
We have $L^G(X,\,\underline\sF,\,t)=1$ for any $\underline\sF\in{\rm NilCoh}_\tau^G(X,\,A)$. The set $1+tA[G][[t]]$ is an abelian group via the multiplication. Then we get a homomorphism of abelian groups
$$L^G(X,\,\bullet,\,t):K_0({\rm Crys}^G(X,\,A))\to 1+tA[G][[t]],\;\;[\underline\sF]\mapsto L^G(X,\,\underline\sF,\,t).$$
In the remainder of this section, we prove the trace formula of equivariant $L$-functions.
\begin{thm}
Let $f:X\to Y$ be a morphism of $k$-schemes. We have a commutative diagram of abelian groups
\[\xymatrix{K_0({\rm Crys}^G(X,\,A))\ar[r]^{Rf_!}\ar[rd]^{L^G(X,\,\bullet,\,t)}&K_0({\rm Crys}^G(Y,\,A))\ar[d]^{L^G(Y,\,\bullet,\,t)}\\& 1+tA[G][[t]].}\]
\end{thm}
\begin{proof}
It suffices to show for any $\underline\sF\in{\rm Coh}^G_\tau(X,\,A)$,
$$L^G(Y,\,Rf_!\underline\sF,\,t)=L^G(X,\,\underline\sF,\,t).$$
By Lemma \ref{zhong}, there is a coherent $\sO_X[G]$-module $\sG$ such that $\sF[G]\simeq\sF\oplus\sG$. Let $\underline\sG=(\sG,\,0)\in{\rm Coh}_\tau^G(X,\,A)$. Then $L^G(X,\,\underline\sF,\,t)=L^G(X,\,\underline\sF\oplus\underline\sG,\,t)$. Let $\tau$ be the composite
$$\bigoplus_{1\leq n\ll\infty}(\sigma_X^n)^*\sF\otimes A[G]=\bigoplus_{1\leq n\ll\infty}(\sigma_X^n)^*\sF[G]\otimes A\x{\tau_\sF\oplus0}\sF[G]\otimes A=\sF\otimes A[G].$$
Then $(\sF,\,\tau)$ is a coherent $\tau$-sheaf over $A[G]$ on $X$. By Lemma \ref{zhong} (2), $$L(X,\,(\sF,\,\tau),\,t)=L^G(X,\,\underline\sF\oplus\underline\sG,\,t)=L^G(X,\,\underline\sF,\,t).$$
So we may assume that $G=\{1\}.$
\end{proof}

(1) Let $0\to\underline\sF'\to\underline\sF\to\underline\sF''\to0$ be a short exact sequence in ${\rm Coh}_\tau(X,\,A)$. If the theorem holds for $\underline\sF'$ and $\underline\sF''$, so does it for $\underline\sF$.

(2) Let $j:U\to X$ be an open immersion with complement $i:Z\to X$. Suppose the theorem holds for $fj$ and $fi$. Then
\begin{eqnarray*}
L(X,\,\underline\sF,\,t)=L(U,\,j^*\underline\sF,\,t)L(Z,\,i^*\underline\sF,\,t)=
L(Y,\,R(fj)_!j^*\underline\sF,\,t)L(Y,\,R(fi)_!i^*\underline\sF,\,t)=L(Y,\,Rf_!\underline\sF,\,t).
\end{eqnarray*}

(3) Suppose $f=gh$ and the theorem holds for $g$ and $h$. Then the theorem holds for $f$.

(4) For any closed point $y\in Y$, consider the cartesian diagram
\[\xymatrix{X_y\ar[r]^{f_y}\ar[d]^{\iota_y}&y\ar[d]^{i_y}\\X\ar[r]^f&Y.}\]
Suppose $L(X_y,\,\iota_y^*\underline\sF,\,t)=L(y,\,R(f_y)_!\iota_y^*\underline\sF,\,t)$ for any $y\in|Y|$. Since $|X|=\coprod_{y\in|Y|}|X_y|$, then by Theorem \ref{yun} (1), we have
\begin{eqnarray*}
L(X,\,\underline\sF,\,t)=\prod_{y\in|Y|}L(X_y,\,\iota_y^*\underline\sF,\,t)=\prod_{y\in|Y|}L(y,\,R(f_y)_!\iota_y^*\underline\sF,\,t)
=\prod_{y\in|Y|}L(y,\,i_y^*Rf_!\underline\sF,\,t)=L(Y,\,Rf_!\underline\sF,\,t).
\end{eqnarray*}

(5) The theorem holds for any finite morphism $f$. In fact, we may assume $Y={\rm Spec}\,k'$ for a finite field extension $k'/k$ and $X$ is a finite $k$-scheme by (4). By (2), we may assume $X$ is connected. Let $i:X_{\rm red}\to X$ be the reduced subscheme of $X$. Then $\underline\sF\simeq i_*i^*\underline\sF\in{\rm Crys}(X,\,A)$. So we may assume that $X$ is reduced. Then $X={\rm Spec}\,k''$ for a finite field extension $k''/k'$. By Definition \ref{3171}, $L({\rm Spec}\,k'',\,\underline\sF,\,t)=L({\rm Spec}\;k',\,f_*\underline\sF,\,t)$.

Let $f$ be a general morphism. By (4) and (2), we may assume that $X$ is affine and $Y={\rm Spec}\,k'$ for some finite field extension $k'/k$. Choose a closed immersion
$i:X\to\A^n_{k'}$. Then $f:X\to{\rm Spec}\,k'$ factors as $X\x{i}\A^n_{k'}\twoheadrightarrow\A^{n-1}_{k'}\twoheadrightarrow\cdots\twoheadrightarrow{\rm Spec}\,k'$.
By (3), we may assume $f$ is the projection $\A^m_{k'}\twoheadrightarrow\A^{m-1}_{k'}$. By (4), we may assume that $f$ is the structure morphism $\A^1_{k''}\to{\rm Spec}\;k''$ for a finite field extension $k''/k'$. We have a cartesian diagram
\[\xymatrix{\A^1_{k''}\ar[r]^\pi\ar[d]^f&\A^1_k\ar[d]^g\\{\rm Spec}\,k''\ar[r]^p&{\rm Spec}\,k.}\]
Apply (5) to the finite morphisms $\pi$ and $p$ and by Theorem \ref{yun} (2), we have $L(\A^1_{k''},\,\underline\sF,\,t)=L(\A^1_k,\,\pi_*\underline\sF,\,t)$ and
$L({\rm Spec \,k''},Rf_!\underline\sF,\,t)=L({\rm Spec}\,k,\,p_*Rf_!\underline\sF,\,t)=L({\rm Spec}\,k,\,Rg_!\pi_*\underline\sF,\,t).$ Then we may assume that $f$ is the structure morphism $\A^1_k\to{\rm Spec}\,k$. Let $\sF_0$ be the torsion part of $\sF$. We have $\tau_\sF\Big(\bigoplus_{1\leq n\ll\infty}(\sigma_X^n)^*\sF_0\otimes A\Big)\subset\sF_0\otimes A$. Then we get a subobject $\underline{\sF_0}$ of $\underline\sF$ supported on a finite closed subscheme $i:Z\to\A^1_k$. Then $fi$ is finite and $\underline{\sF_0}=i_*i^*\underline{\sF_0}$. By (5), we have $$L(X,\,\underline{\sF_0},\,t)=L(Z,\,i^*\underline{\sF_0},\,t)=L({\rm Spec}\,k,(fi)_*i^*\underline{\sF_0},\,t)=L({\rm Spec}\,k,Rf_!i_*i^*\underline{\sF_0},t)=L({\rm Spec}\,k,\,Rf_!\underline{\sF_0},\,t).$$
Applying (1) to the short exact sequence $0\to\underline{\sF_0}\to\underline\sF\to\underline\sF/\underline{\sF_0}\to0$, we may assume that $\sF$ is torsion free over $\sO_{\A^1_k}$ and then $\sF\simeq\sO_{\A^1_k}^{\oplus r}$ for some $r$. Let $j:\A^1_k\hookrightarrow\P^1_k$ be the inclusion and $\bar f$ the structure morphism $\P^1_k\to{\rm Spec}\,k$. Then $\tau_{j_*\sF}\Big(\bigoplus_{1\leq n\ll\infty}(\sigma^n_{\P^1_k})^*\sO^{\oplus r}_{\P^1_k}(-d\infty)\otimes A\Big)\subset\sO_{\P^1_k}^{\oplus r}(-d\infty)\otimes A$ for $d\gg0$. Then $j_!\underline\sF$ is represented by $(\sO^r_{\P^1_k}(-d\infty),\,\tau_{\sO^{\oplus r}_{\P_k^1}(-d\infty)})$ for $d\gg0$, where $\tau_{\sO^{\oplus r}_{\P_k^1}(-d\infty)}:\bigoplus_{1\leq n\ll\infty}(\sigma^n_{\P^1_k})^*\sO^{\oplus r}_{\P^1_k}(-d\infty)\otimes A\to\sO_{\P^1_k}^{\oplus r}(-d\infty)\otimes A$ is induced by $\tau_{j_*\sF}$.
For $d\gg0$, $$R^0\bar f_*\sO^{\oplus r}_{\P_k^1}(-d\infty)=0\hbox{ and }R^1\bar f_*j_!\underline\sF=\Big(H^1(\P^1_k,\,\sO_{\P^1_k}^{\oplus r}(-d\infty)),\tau_{R^1\bar f_*\sO^{\oplus r}_{\P_k^1}(-d\infty)}\Big).$$
By Serre duality, we have
$$H^1(\P^1_k,\,\sO_{\P^1_k}^{\oplus r}(-d\infty))={\rm Hom}_k\Big(H^0(\P^1_k,\,\Omega_{\P^1_k/k}^1(d\infty))^{\oplus r},\,k\Big).$$
Then $\tau_{R^1\bar f_*\sO^{\oplus r}_{\P_k^1}(-d\infty)}$ induces a family of operators $\{C_n\}_{n\geq1}$
on $H^0(\P^1_k,\,\Omega_{\P^1_k/k}^1(d\infty))^{\oplus r}\otimes A$.
By $$\varinjlim_d H^0(\P^1_k,\,\Omega_{\P^1_k/k}^1(d\infty))^{\oplus r}={\rm Hom}_{\sO_{\A^1_k}}(\sO_{\A^1_k}^{\oplus r},\,\Omega_{\A^1_k/k}^1),$$ $C_n$ induces an $n$-th Cartier operator on ${\rm Hom}_{\sO_{\A^1_k}}(\sO_{\A^1_k}^{\oplus r},\,\Omega_{\A^1_k/k}^1)\otimes A$ which coincides with the Cartier operator defined in Theorem \ref{ji}.
Moreover, $H^0(\P^1_k,\,\Omega_{\P^1_k/k}^1(d\infty))^{\oplus r}$ is a common nucleus for $\{C_n\}_{n\geq1}$ on ${\rm Hom}_{\sO_{\A^1_k}}(\sO_{\A^1_k}^{\oplus r},\,\Omega^1_{\A^1_k/k})\otimes A$ for $d\gg0$. By Theorem \ref{ji},
we have
\begin{eqnarray*}
L(\A^1_k,\,\underline\sF,\,t)&=&{\rm det}_{A[t]}\Big(1-\sum_{1\leq n\ll\infty}t^nC_n,\,H^0(\P^1_k,\,\Omega^1_{\P^1_k/k}(d\infty))^{\oplus r}\otimes A[t]\Big)\\
&=&{\rm det}_{A[t]}\Big(1-\sum_{1\leq n\ll\infty}t^n\tau_n,\,H^1(\P_k^1,\,\sO_{\P^1_k}^{\oplus r}(-d\infty))\otimes A[t]\Big)\\
&=&L({\rm Spec}\,k,\,Rf_!\underline\sF,\,t).
\end{eqnarray*}
This completes the proof.
\subsection{Another type of $L$-functions}
\begin{defn}
(1) A coherent $\widetilde\tau$-sheaf over $A$ on $X$ is a pair $(\sF,\,\widetilde\tau_\sF)$ consisting of a coherent $\sO_X$-module $\sF$ and an $\sO_X\otimes A$-linear homomorphism
$$\widetilde\tau_\sF:\bigoplus_{0\leq n\ll\infty}(\sigma_X^n)^*\sF\otimes A\to\sF\otimes A.$$
Then a coherent $\widetilde\tau$-sheaf $(\sF,\,\widetilde\tau_\sF)$ over $A$ on $X$ is a coherent $\tau$-sheaf $(\sF,\,\tau_\sF)$ over $A$ on $X$ together with $\tau_0\in{\rm End}_{\sO_X\otimes A}(\sF\otimes A)$, where $\tau_\sF:\bigoplus_{1\leq n\ll\infty}(\sigma_X^n)^*\sF\otimes A\to\sF\otimes A$ is truncated by $\widetilde\tau_\sF$.

(2) A coherent $\widetilde\tau$-sheaf over $A$ on $X$ is called nilpotent if the associated $\tau$-sheaf is nilpotent.

(3) Let ${\rm Coh}_{\widetilde\tau}^G(X,\,A)$ (resp. ${\rm NilCoh}^G_{\widetilde\tau}(X,\,A)$) be the category of $G$-modules in the category of coherent (resp. nilpotent) $\widetilde\tau$-sheaves over $A$ on $X$. Let $\widetilde{\rm Crys}^G(X,\,A)=\frac{{\rm Coh}_{\widetilde\tau}^G(X,\,A)}{{\rm NilCoh}^G_{\widetilde\tau}(X,\,A)}$.
\end{defn}
\begin{defn}\label{3163}
Suppose $G$ is a finite abelian group of order prime to $p$.

(1) Let $x={\rm spec}\;k'$ for a finite field extension $k'/k$. Any $(\sF,\,\widetilde\tau_\sF)\in{\rm Coh}_{\widetilde\tau}^G(x,\,A)$ consists of an $n$-th Frobenius operator $\tau_n$ on $\sF\otimes_kA$ over $A$ on $k'$ for finitely many $n\in\N$. Define
$$L^G(x,\,(\sF,\widetilde\tau_\sF),\,t)=\frac{{\rm det}_{A[G][[t^{-1}]]}\Big(1-t^{-1}\tau_0,\;\sF\otimes_kA[[t^{-1}]]\Big)}{{\rm det}_{A[G][[t^{-1}]]}\Big(1-t^{-1}\sum_{0\leq n\ll\infty}\tau_n,\;\sF\otimes_kA[[t^{-1}]]\Big)}\in 1+t^{-1}A[G][[t^{-1}]].$$
(2) For any $(\sF,\,\widetilde\tau_\sF)\in{\rm Coh}_{\widetilde\tau}^G(X,\;A)$, define
$$L^G(X,\,(\sF,\widetilde\tau_\sF),\,t)=\prod_{x\in|X|}L^G(x,\,i_x^*(\sF,\,\widetilde\tau_\sF),\,t).$$
\end{defn}
\begin{lem}\label{3172}
Suppose $G$ is an abelian group of order prime to $p$. We get a homomorphism of abelian groups
$$L^G(X,\,\bullet,\,t):K_0(\widetilde{\rm Crys}^G(X,\,A))\to 1+t^{-1}A[G][[t^{-1}]],\;\;[(\sF,\,\widetilde\tau_\sF)]\mapsto L^G(X,\,(\sF,\,\widetilde\tau_\sF),\,t).$$
\end{lem}
\section{$v$-adic $L$-values of shtukas on curves}
\noindent In \cite{L}, Lafforgue studied the $v$-adic $L$-values of shtukas on curves. In this section, we generalize it to an equivariant version. Let $G$ be a finite abelian group of order prime to $p$. Let $X$ be a smooth projective curve and $T={\rm Spec}\;A$ a smooth affine curve over $k$. Let $i:\sE\to\sE'$ be a morphism of $G$-bundles on $X\times T$ such that $i$ is isomorphic at the generic point. Let $r$ be a positive integer. For any $1\leq s\leq r$, let $\tau_s:\sE\to\sE'$ be an $s$-th Frobenius map over $A$ on $X$.

Let $Z(\det(i))$ be the zeros of $\det(i)$ in $X\times T$. Fix $v\in|T|$ such that $Z(\det(i))\cap X\times\{v\}$ is finite. Take a finite subset $S$ of $|X|$ such that $S\times\{v\}\supset Z(\det(i))\cap X\times\{v\}$. Let $A_v$ be the completion of $A$ at $v$ and choose a uniformizer of $A_v$ which is also denoted by $v$.
For any $x\in|X|-S$, let $i_x,\,(\tau_1)_x,\ldots,(\tau_r)_x:\sE_x\to\sE'_x$ be the restriction of $i,\,\tau_1,\ldots,\tau_r$
on ${\rm Spec}\;k(x)\otimes_k A_v$, respectively.
Then $i_x:\sE_x\to\sE_x'$ is isomorphic for any $x\in|X|-S$.

\begin{lem}
Define the $v$-adic $G$-equivariant $L$-function away from $S$ of the diagram
$\sE\x{i}\sE'\xleftarrow{j=(\tau_1,\ldots,\tau_r)}\sE$
to be
$$L_{v}^G(X-S,\,(\sE,\sE',i,j),T)=\prod_{x\in|X|-S}{\rm det}_{A_v[G]}\Big(1-\sum_{s=1}^rT^si_x^{-1}(\tau_s)_x,\,\sE_x\Big)^{-1}\in 1+TA_v[G][[T]].$$
Then $L_{v}^G(X-S,\,(\sE,\sE',i,j),T)\in 1+TA_v[G]\langle\langle T\rangle\rangle$.
\end{lem}
\begin{proof}

For any $n\geq1$, let $i_{n,\,x},\,(\tau_1)_{n,\,x},\ldots,(\tau_r)_{n,\,x}:\sE_{n,\,x}\to\sE'_{n,\,x}$ be the restriction of $i,\,\tau_1,\ldots,\tau_r$
on ${\rm Spec}\;k(x)\otimes_k A/v^n$, respectively. If $x\in |X|-S$, then $i_{n,\,x}$ is an isomorphism. By Definition \ref{yuan} and Theorem \ref{ji}, we have
$$\prod_{x\in|X|-S}{\rm det}_{A/v^{n}[G]}\Big(1-\sum_{s=1}^rT^si_{n,\,x}^{-1}(\tau_s)_{n,\,x},\,\sE_{n,\,x}\Big)^{-1}\in 1+TA/v^n[G][T].$$
This proves the lemma.
\end{proof}

\begin{defn}
We have
$$L_{v}^G(X-S,\,(\sE,\sE',i,j),T)=(1-T)^sg(T)$$ for some $s\in\N$ and $g(T)\in 1+TA_v[G]\langle\langle T\rangle\rangle$ such that
$0\neq g(1)\in A_v[G]$. The $G$-equivariant $v$-adic $L$-value $L_{v}^{G}(X-S,\,(\sE,\sE',i,j))$
of the diagram $(\sE,\sE',i,j)$ away from $S$ is defined to be $g(1).$
\end{defn}

Suppose $X-S={\rm Spec}\;R_S$. Let $\sM$ and $\sM'$ be the $R_S\widehat\otimes A_v$-modules defined by $\sE$ and $\sE'$. For any $w\in S$, let $\sO_w$ be the complete local ring of $X$ at $w$. Let $\sM_w$ and $\sM_w'$ be the $\sO_w\widehat\otimes A_v$-modules defined by $\sE$ and $\sE'$. Let $\sV={\rm Hom}_{R_S\widehat\otimes A_v}(\sM,\Omega_{R_S/k}^1\widehat\otimes A_v)$ and $\sV'={\rm Hom}_{R_S\widehat\otimes A_v}(\sM',\Omega_{R_S/k}^1\widehat\otimes A_v)$. Let $\widetilde j=\sum_{s=1}^r\tau_s$.
We have two commutative diagrams in the derived category of $A_v[G]$-modules
\begin{equation}\label{com}
 \xymatrix{R\Gamma(X,\,\sE)\otimes_AA_v\ar[r]\ar@<0.5ex>[d]_{i~~}\ar@<-0.5ex>[d]
^{~~i-\widetilde j}&\bigoplus_{w\in S}\sM_w\ar[r]\ar@<0.5ex>[d]_{i~~}\ar@<-0.5ex>[d]^{~~i-\widetilde j}&{\rm Hom}_{A_v}(\sV,\,A_v)\ar@<0.5ex>[d]_{i~~}\ar@<-0.5ex>[d]^{~~i-\widetilde j}\\
R\Gamma(X,\,\sE')\otimes_AA_v\ar[r]&\bigoplus_{w\in S}\sM_w'\ar[r]&{\rm Hom}_{A_v}(\sV',\,A_v),}
\end{equation}
one with the left arrows and one with the right arrows.
Here the morphism $\bigoplus_{w\in S}\sM_w\to{\rm Hom}_{A_v}(\sV,\,A_v)$ (resp. $\bigoplus_{w\in S}\sM'_w\to{\rm Hom}_{A_v}(\sV',\,A_v)$)
associates each $(f_w)\in\bigoplus\sM_w$ and $g\in\sV$ (resp. $(f_w)\in\bigoplus\sM_w'$ and $g\in\sV'$) to the sum of residue of $\langle g,\,f_w\rangle$ at $w$.
\begin{lem}\label{j}
Suppose for any $w\in S$, there exists two $A_v[G]$-linear isomorphisms
$\exp:\sM_w\simeq\sM_w$ and $\exp:\sM_w'\simeq\sM_w'$ which satisfy the following three conditions.
\begin{enumerate}
\item  $\exp\circ i=i-\widetilde j\circ\exp:\sM_w\to\sM_w'$.
\item For any $t\in\N$, $\exp(w^t\sM_w)=w^t\sM_w$, $\exp(w^t\sM_w')=w^t\sM_w'$, $(\exp-{\rm id})(w^t\sM_w)\subset w^{t+1}\sM_w$ and $(\exp-{\rm id})(w^t\sM'_w)\subset w^{t+1}\sM'_w$.
\item For any $s\in\N$,  $(\exp-{\rm id})(w^t\sM_w)\subset w^{t+s}\sM_w$ and $(\exp-{\rm id})(w^t\sM'_w)\subset w^{t+s}\sM'_w$ for $t$ large enough.
\end{enumerate}
Let $\log:\sM_w\simeq\sM_w$ and $\log:\sM_w'\simeq\sM_w'$ be the inverse maps of $\exp$. Let $\bigoplus_{w\in S}\sM_w'\x{\pi}C$ be the cokernel of $\bigoplus_{w\in S}\sM_w\x{i}\bigoplus_{w\in S}\sM_w'$. Denote by $\iota$ the natural map $R\Gamma(X,\,\sE')\otimes_AA_v\to\bigoplus_{w\in S}\sM_w'$.
Suppose $i-\widetilde j:{\rm Hom}_{A_v}(\sV,\,A_v)\to{\rm Hom}_{A_v}(\sV',\,A_v)$ is an isomorphism. Then we have a commutative diagram
\begin{eqnarray}\label{3202}\xymatrix{R\Gamma(X,\,\sE)\otimes_AA_v\ar[r]\ar@<0.5ex>[d]_{i~~}\ar@<-0.5ex>[d]
^{~~i-\widetilde j}&\bigoplus_{w\in S}\sM_w\ar[r]\ar@<0.5ex>[d]_{i~~}\ar@<-0.5ex>[d]^{~~i-\widetilde j}&{\rm Hom}_{A_v}(\sV,\,A_v)\ar@<0.5ex>[d]_{i~~}\ar@<-0.5ex>[d]^{~~i-\widetilde j}\\
R\Gamma(X,\,\sE')\otimes_AA_v\ar[r]^\iota\ar@<0.5ex>[d]_{\pi\iota~~}\ar@<-0.5ex>[d]^{~~\pi\log\iota}&\bigoplus_{w\in S}\sM_w'\ar[r]\ar@<0.5ex>[d]_{\pi~~}\ar@<-0.5ex>[d]^{~~\pi\log}
&{\rm Hom}_{A_v}(\sV',\,A_v)\\
C\ar@{=}[r]&C},\end{eqnarray}
whose four vertical triangles are distinguished.
Let $\delta$ be the composite
\begin{eqnarray*}
{\rm det}_{A_v[G]}(C)\simeq{\rm det}_{A_v[G]}\Big(R\Gamma(X,\,\sE)\otimes_AA_v\Big)^{-1}\bigotimes{\rm det}_{A_v[G]}\Big(R\Gamma(X,\,\sE')\otimes_AA_v\Big)\simeq{\rm det}_{A_v[G]}(C),
\end{eqnarray*}
where the first and the second isomorphism are given by the first and the second vertical distinguished triangles of diagram (\ref{3202}), respectively.
We have
$$L_{v}^G(X-S,\,(\sE,\sE',i,j))=[\delta:1]\in A_v[G]^\times.$$
\end{lem}
\begin{proof}
This lemma is an equivariant version of \cite[ Lemma 2.5 ] {F} and the proof is highly similar to that of \cite[ Lemma 2.5 ] {F} by using the main results established in section 2 and section 3.
\end{proof}
\section{Proofs of Theorem \ref{391} and Theorem \ref{18}}
Recall that $K$ is a finite extension of $k(t)$ and $L/K$ is a finite Galois extension of Galois group $G$. Let $K_\infty=K\otimes_{k(t)}k((t^{-1}))$ and $L_\infty=L\otimes_{k(t)}k((t^{-1}))$. Let $X$ and $\widetilde X$ be the projective curves over $k$ of function fields $K$ and $L$, respectively. Let $\sO_\infty=\sO_X\otimes_{\sO_{\P^1_k}}k[[t^{-1}]]$ and $\widetilde\sO_\infty=\sO_{\widetilde X}\otimes_{\sO_{\P_k^1}}k[[t^{-1}]]$.
Let $E$ be an abelian $t$-module over $\sO_K$ of dimension $n$ defined by
$$E(t)=A_0+A_1\tau+\ldots A_r\tau^r\in M_n(\sO_K)\{\tau\}$$
such that $(A_0-zI_n)^n=0$ where $z$ is the image of $t$ in $\sO_K$.

In this section, we mainly study the special $L$-value $L(E,\,\rho)$ for a finite dimensional $F$-linear representation $\rho:G\to{\rm GL}_F(V)$. If $G$ is an abelian group of order prime to $p$, we also study the equivariant special $L$-value $L(E,\,G)$.

For any coherent sheaf $\sF$ on a $\P_k^1$-scheme $Z$ and any $d\in\Z$, let $\sF(d\infty)=\sF\otimes_{\sO_Z}\sO_{\P_k^1}(d\infty)|_Z$.
View $\sO_{\widetilde X}$ as a vector bundle on $X$. Since $\sO_\infty$ is flat over $\sO_X$ and $K_\infty$ is flat over $\sO_K$, we have
\begin{eqnarray*}
&&{\rm Hom}_{k[G]}(V,\,\sO_{\widetilde X}^n(-d\infty))\otimes_{\sO_X}\sO_\infty\simeq{\rm Hom}_{k[G]}(V,\,\sO_{\widetilde X}^n(-d\infty)\otimes_{\sO_X}\sO_\infty)\simeq{\rm Hom}_{k[G]}(V,\,z^{-d}\widetilde\sO^n_\infty);\\
&&{\rm Hom}_{k[G]}(V,\,\sO_L^n)\otimes_{\sO_K}K_\infty={\rm Hom}_{k[G]}(V,\,\sO_L^n\otimes_{\sO_K}K_\infty)={\rm Hom}_{k[G]}(V,\,L_\infty^n).
\end{eqnarray*}
Since ${\rm Hom}_{k[G]}(V,\,\sO_L^n)$ is a flat $\sO_K$-module, we have
$${\rm Hom}_{k[G]}(V,\,\sO_L^n)\otimes_{\sO_K}\frac{K_\infty}{\sO_K}=\frac{{\rm Hom}_{k[G]}(V,\,L_\infty^n)}{{\rm Hom}_{k[G]}(V,\,\sO_L^n)}.$$

Let $\sE_\rho={\rm Hom}_{k[G]}(V,\,\sO_{\widetilde X}^n(-d\infty))$ and $\sE_G=\sO_{\widetilde X}^n(-d\infty)$.
Then $\sE_\rho$ is a vector bundle on $X$ and $\sE_G$ is a $G$-bundle on $X$. For any $k$-algebra $R$, let $R_\rho=F\otimes_kR$ and $R_G=R[G]$. For any $k[G]$-module $M$, let $M_\rho={\rm Hom}_{k[G]}(V,\,M)$ and $M_G=M$.

Let $*=\rho$ or $G$. Then we have $\sE_*\otimes_{\sO_X}\sO_\infty=(z^{-d}\widetilde\sO_\infty^n)_*$ and
$$\Gamma({\rm Spec}\,\sO_K,\,\sE_*)\otimes_{\sO_K}\frac{K_\infty}{\sO_K}=(\sO_L^n)_*\otimes_{\sO_K}\frac{K_\infty}{\sO_K}=\frac{(L_\infty^n)_*}{(\sO_L^n)_*}.$$
For any vector bundle $\sE$ on $X$, we have a distinguished triangle of complexes of $k$-modules
\begin{eqnarray*}
R\Gamma(X,\,\sE)\to\sE\otimes_{\sO_X}\sO_\infty\to\Gamma({\rm Spec}\,\sO_K,\,\sE)\otimes_{\sO_K}\frac{K_\infty}{\sO_K}.
\end{eqnarray*}
In particular, for $*=\rho$ or $G$, we have a distinguished triangle
$$R\Gamma(X,\,\sE_*)\to(z^{-d}\widetilde\sO_\infty^n)_*\to\frac{(L_\infty^n)_*}{(\sO_L^n)_*}.$$

Choose $e\in\N$ such that $A_0\in M_n(\Gamma(X,\sO_X(e\infty)))$. Then $e\geq1$.
We get two $k[G]$-linear maps $1-t^{-1}A_0$ and $1-t^{-1}\sum_{s=0}^rA_s\tau^s:\sO^n_{\widetilde X}(-d\infty)((t^{-1}))\to\sO^n_{\widetilde X}((e-d)\infty)((t^{-1}))$ and $z^{-d}\widetilde \sO_\infty^n((t^{-1}))\to z^{e-d}\widetilde\sO_\infty^n((t^{-1}))$. By \cite [section 3] {F}, $\exp_E$ defines a $k[G]$-linear automorphism of $z^{-d}\widetilde\sO_\infty^n$ and $z^{e-d}\widetilde\sO_\infty^n$ for $d\gg0$.
It induces a $k((t^{-1}))_*$-linear automorphism on $(z^{-d}\widetilde\sO_\infty^n)_*((t^{-1}))$ and $(z^{e-d}\widetilde\sO_\infty^n)_*((t^{-1}))$ for $d\gg0$ which is also denoted by $\exp_E$. For $d\gg0$, we have
$$\exp_E\circ(1-t^{-1}A_0)=(1-t^{-1}\sum_{s=0}^rA_s\tau^s)\circ\exp_E:(z^{-d}\widetilde\sO_\infty^n)_*((t^{-1}))\to (z^{e-d}\widetilde\sO_\infty^n)_*((t^{-1})).$$
For any $s\in\N$, $(\exp_E-{\rm id})((z^{-d}\widetilde\sO_\infty^n)_*)\subset (z^{-d-s}\widetilde\sO_\infty^n)_*$ for $d\gg0$.
Let $\log_E$ be the inverse map of $\exp_E$ on $(z^{-d}\widetilde\sO_\infty^n)_*((t^{-1}))$ and $(z^{e-d}\widetilde\sO_\infty^n)_*((t^{-1}))$.

Let $(z^{e-d}\widetilde\sO^n_\infty)_*((t^{-1}))\x{\pi}C$ be the cokernel of $(z^{-d}\widetilde\sO^n_\infty)_*((t^{-1}))\x{1-t^{-1}A_0}(z^{e-d}\widetilde\sO^n_\infty)_*((t^{-1}))$.
Then $(z^{e-d}\widetilde\sO^n_\infty)_*((t^{-1}))\x{\pi\log_E}C$ is the cokernel of $(z^{-d}\widetilde\sO^n_\infty)_*((t^{-1}))\x{1-t^{-1}\sum_{s=0}^rA_s\tau^s}(z^{e-d}\widetilde\sO^n_\infty)_*((t^{-1})).$
Since $1-t^{-1}A_0$ and $1-t^{-1}\sum_{s=0}^rA_s\tau^s$ on $\frac{(L_\infty^n)_*}{(\sO_L^n)_*}((t^{-1}))$ are isomorphic, we have a commutative diagram
\begin{eqnarray}\label{3191}
\xymatrix{R\Gamma(X,\,\sE_*)((t^{-1}))\ar[r]\ar@<0.5ex>[d]_{1-t^{-1}A_0~~}\ar@<-0.5ex>[d]
^{~~1-t^{-1}\sum_{s=0}^rA_s\tau^s }&(z^{-d}\widetilde\sO_\infty^n)_*((t^{-1}))\ar[rr]\ar@<0.5ex>[d]_{1-t^{-1}A_0~~}\ar@<-0.5ex>[d]^{~~1-t^{-1}\sum_{s=0}^rA_s\tau^s}
&&\frac{(L_\infty^n)_*}{(\sO_L^n)_*}((t^{-1}))
\ar@<0.5ex>[d]_{1-t^{-1}A_0~}\ar@<-0.5ex>[d]^{~1-t^{-1}\sum_{s=0}^rA_s\tau^s}\\
R\Gamma(X,\,\sE_*(e\infty))((t^{-1}))\ar[r]^\iota\ar@<0.5ex>[d]_{\pi\iota~~}\ar@<-0.5ex>[d]^{~~\pi\log_E\iota}&(z^{e-d}\widetilde\sO_\infty^n)_*((t^{-1}))
\ar[rr]\ar@<0.5ex>[d]_{\pi~~}\ar@<-0.5ex>[d]^{~~\pi\log_E}
&&\frac{(L_\infty^n)_*}{(\sO_L^n)_*}((t^{-1}))\\
C\ar@{=}[r]&C,}
\end{eqnarray}
whose four vertical triangles are distinguished.

\begin{lem}\label{c1}
For $d\gg0$, we have two distinguished triangles
\begin{eqnarray}\label{f1}
&&R\Gamma(X,\,\sE_*)((t^{-1}))\x{t-A_0}R\Gamma\big(X,\,\sE_*(e\infty)\big)((t^{-1}))\x{\pi\iota}C;\\\label{f2}
&&R\Gamma(X,\,\sE_*)((t^{-1}))\x{t-\sum_{s=0}^rA_s\tau^s}R\Gamma\big(X,\,\sE_*(e\infty)\big)((t^{-1}))\x{\pi\log_E\iota}C.
\end{eqnarray}
Then $L(E,\;*)=[\delta_1:1]$ where $\delta_1$ is the composite
\begin{eqnarray*}
{\rm det}_{k((t^{-1}))_*}(C)&\simeq&{\rm det}_{k((t^{-1}))_*}\Big(R\Gamma(X,\,\sE_*)((t^{-1}))\Big)^{-1}\bigotimes{\rm det}_{k((t^{-1}))_*}\Big(R\Gamma\big(X,\,\sE_*(e\infty)\big)((t^{-1}))\Big)\\&\simeq&{\rm det}_{k((t^{-1}))_*}(C)
\end{eqnarray*}
where the first and the second isomorphisms are defined by (\ref{f1}) and (\ref{f2}), respectively.
\end{lem}
\begin{proof}
For $d\gg0$, we easily get $[\delta_1:1]=[\delta_2:1]\in k((t^{-1}))_*$, where $\delta_2$ is the composite
\begin{eqnarray*}
{\rm det}_{k((t^{-1}))_*}(C)&\simeq&
{\rm det}_{k((t^{-1}))_*}\Big(R\Gamma(X,\,\sE_*)((t^{-1}))\Big)^{-1}\bigotimes{\rm det}_{k((t^{-1}))_*}\Big(R\Gamma\big(X,\,\sE_*(e\infty)\big)((t^{-1}))\Big)
\\&\simeq&{\rm det}_{k((t^{-1}))_*}(C),
\end{eqnarray*}
where the first and the second isomorphisms are defined by the first and the second vertical triangles of diagram (\ref{3191}), respectively.
Applying Theorem \ref{3155} to the $\widetilde\tau$-sheaf $({\rm Hom}_k(V,\,\sO_L^n),\;A_0,A_1\tau,\ldots,A_r\tau^r)\in{\rm Coh}_{\widetilde\tau}^G({\rm Spec}\;\sO_K,\,F)$, we have
\begin{eqnarray}\label{3156}
{\rm Hom}_{k[G]}(V,\,\sO_L^n)\otimes_{\sO_K}\sO_K/\mathfrak p\simeq{\rm Hom}_{k[G]}(V,\,\sO_L^n/\mathfrak p\sO_L^n)\in\widetilde{\rm Crys}({\rm Spec}\;\kappa_\mathfrak p,\,F).
\end{eqnarray}
%By Lemma \ref{3172}, we have
%\begin{eqnarray*}
%&&\frac{{\rm det}_{F[t]}\Big(t-A_0,\;{\rm Hom}_{k[G]}(V,\,\sO_L)^n\otimes_{\sO_K}\sO_K/\mathfrak p[t]\Big)}{{\rm det}_{F[t]}
%\Big(t-\sum_{s=0}^rA_s\tau^s,\;{\rm Hom}_{k[G]}(V,\,\sO_L)^n\otimes_{\sO_K}\sO_K/\mathfrak p[t]\Big)}\\
%&=&\frac{{\rm det}_{F[t]}\Big(t-A_0,\;{\rm Hom}_{k[G]}(V,\,\sO_L^n/\mathfrak p\sO_L^n)[t]\Big)}{{\rm det}_{F[t]}
%\Big(t-\sum_{s=0}^rA_s\tau^s,\;{\rm Hom}_{k[G]}(V,\,\sO_L^n/\mathfrak p\sO_L^n)[t]\Big)}.
%\end{eqnarray*}
For $d\gg0$, the diagram
$$\sE_*[t^{-1}]\x{1-t^{-1}A_0}\sE_*(e\infty)[t^{-1}]\xleftarrow{t^{-1}A_1\tau,\ldots,t^{-1}A_r\tau^r}\sE_*[t^{-1}]$$
satisfies the assumptions of Lemma \ref{j} at the prime ideal of $k[t^{-1}]$ generated by $t^{-1}$. By Lemma \ref{j}, we have
\begin{eqnarray*}
&&[\delta_2:1]\\&=&\prod_{\mathfrak p\in{\rm Max}(\sO_K)}{\rm det}_{k[[t^{-1}]]_*}\Big(1-(1-t^{-1}A_0)^{-1}\sum_{s=1}^rt^{-1}A_s\tau^s,\;(\sO_L^n)_*\otimes_{\sO_K}\sO_K/\mathfrak p[[t^{-1}]]\Big)^{-1}\\
&=&\prod_{\mathfrak p\in{\rm Max}(\sO_K)}\frac{{\rm det}_{k[[t^{-1}]]_*}\Big(1-t^{-1}A_0,\;(\sO_L^n)_*\otimes_{\sO_K}\sO_K/\mathfrak p[[t^{-1}]]\Big)}{{\rm det}_{k[[t^{-1}]]_*}\Big(1-t^{-1}\sum_{s=0}^rA_s\tau^s,\;(\sO_L^n)_*\otimes_{\sO_K}\sO_K/\mathfrak p[[t^{-1}]]\Big)}\\
&=&\prod_{\mathfrak p\in{\rm Max}(\sO_K)}\frac{{\rm det}_{k[[t^{-1}]]_*}\Big(1-t^{-1}A_0,\;(\sO_L^n/\mathfrak p\sO_L^n)_*[[t^{-1}]]\Big)}{{\rm det}_{k[[t^{-1}]]_*}
\Big(1-t^{-1}\sum_{s=0}^rA_s\tau^s,\;(\sO_L^n/\mathfrak p\sO_L^n)_*[[t^{-1}]]\Big)}\\
&=&\prod_{\mathfrak p\in{\rm Max}(\sO_K)}\frac{{\rm det}_{k[t]_*}\Big(t-A_0,\;(\sO_L^n/\mathfrak p\sO_L^n)_*[t]\Big)}{{\rm det}_{k[t]_*}
\Big(t-\sum_{s=0}^rA_s\tau^s,\;(\sO_L^n/\mathfrak p\sO_L^n)_*[t]\Big)}\\
&=&\prod_{\mathfrak p\in{\rm Max}(\sO_K)}\frac{|{\rm Lie}(E)(\sO_L/\mathfrak p\sO_L)_*|_{k[t]_*}}{|E(\sO_L/\mathfrak p\sO_L)_*|_{k[t]_*}}\\
&=&L(E,\,*)\in 1+t^{-1}k[[t^{-1}]]_*,
\end{eqnarray*}
where the third equality holds by applying Lemma \ref{3172} to (\ref{3156}).
\end{proof}

By \cite [Lemma 1.7 ]{F}, the surjective $k[t]$-linear map
$$z^{e-d}\widetilde\sO_\infty^n((t^{-1}))\to{\rm Lie}(E)(L_\infty),\;\;\sum_{s\ll+\infty}x_st^s\mapsto\sum_{s\ll+\infty}A_0^s(x_s)$$
induces a $k((t^{-1}))_*$-linear map
$$q:(z^{e-d}\widetilde\sO_\infty^n)_*((t^{-1}))\to{\rm Lie}(E)(L_\infty)_*.$$
Let $\widetilde f\in{\rm Hom}_{k[G]}(V,\,{\rm Lie}(E)(L_\infty))$. Since $V$ is a finite set, then $z^{-sq^n}\widetilde f(V)\subset z^{e-d}\widetilde\sO^n_\infty$ for $s\gg0$.
Define $f:V\to z^{e-d}\widetilde\sO^n_\infty((t^{-1})),\;v\mapsto z^{-sq^n}\widetilde f(v)t^{sq^n}.$ The fact $A_0^{q^n}=z^{q^n}$ shows $q(f)=\widetilde f$ and hence $q$ is surjective for $*=\rho$. Actually, the composite $(z^{e-d}\widetilde\sO_\infty^n)_*[t]\hookrightarrow(z^{e-d}\widetilde\sO_\infty^n)_*((t^{-1}))\x{q}{\rm Lie}(E)(L_\infty)_*$ is also surjective.
By the exact sequence
$$(z^{-d}\widetilde\sO_\infty^n)_*((t^{-1}))\x{1-t^{-1}A_0}(z^{e-d}\widetilde\sO_\infty^n)_*((t^{-1}))\x{\pi}C\to0$$ and $q\circ (1-t^{-1}A_0)=0$,
$q$ factors as $(z^{e-d}\widetilde\sO_\infty^n)_*((t^{-1}))\x{\pi}C\x{p}{\rm Lie}(E)(L_\infty)_*$.
Let $(z^{e-d}\widetilde\sO_\infty^n)_*[t]\x{\pi_1}D$ be the cokernel of $(z^{-d}\widetilde\sO_\infty^n)_*[t]\x{t-A_0}(z^{e-d}\widetilde\sO_\infty^n)_*[t].$
For $d\gg0$, automorphisms $\log_E$ on $z^{-d}\widetilde\sO_\infty^n[t]$ and $z^{e-d}\widetilde\sO_\infty^n[t]$ induce
$$\log_E\circ(t-\sum_{s=0}^rA_s\tau^s)=(t-A_0)\circ\log_E:(z^{-d}\widetilde\sO_\infty^n)_*[t]\to (z^{e-d}\widetilde\sO_\infty^n)_*[t].$$
Then $(z^{e-d}\widetilde\sO_\infty^n)_*[t]\x{\pi_1\log_E}D$ is the cokernel of $(z^{-d}\widetilde\sO_\infty^n)_*[t]\x{t-\sum_{s=0}^rA_s\tau^s}(z^{e-d}\widetilde\sO_\infty^n)_*[t].$
There exists a unique $k[t]_*$-linear map $\eta:D\to C$ such that the following diagram commutes:
\[\xymatrix{(z^{-d}\widetilde\sO_\infty^n)_*[t]\ar[r]^{t-A_0}\ar[d]&(z^{e-d}\widetilde\sO_\infty^n)_*[t]\ar[rr]^{\pi_1}
\ar[d]&&D\ar[d]^\eta\ar[rr]^{p\eta}&&{\rm Lie}(E)(L_\infty)_*\ar@{=}[d]\\
(z^{-d}\widetilde\sO_\infty^n)_*((t^{-1}))\ar[r]^{t-A_0}&(z^{e-d}\widetilde\sO_\infty^n)_*((t^{-1}))\ar[rr]^\pi &&C\ar[rr]^p&&{\rm Lie}(E)(L_\infty)_*.}\]
We have shown that $p\eta\pi_1$ is surjective and so are $p\eta$ and $p$.
By \cite[ Lemma 3.2 ]{F}, we have two short exact sequences
\begin{eqnarray*}\label{10}
&&0\to\frac{(L_\infty^n)_*}{(\sO_L^n)_*}[t]\x{t-A_0}\frac{(L_\infty^n)_*}{(\sO_L^n)_*}[t]\x{}\frac{{\rm Lie}(E)(L_\infty)_*}{{\rm Lie}(E)(\sO_L)_*}\to0;\\\label{11}
&&0\to\frac{(L_\infty^n)_*}{(\sO_L^n)_*}[t]\x{t-\sum_{s=0}^rA_s\tau^s}\frac{(L_\infty^n)_*}{(\sO_L^n)_*}[t]\x{}
\frac{E(L_\infty)_*}{E(\sO_L)_*}\to0
\end{eqnarray*}
Let $\bar p$ be the natural map ${\rm Lie}(E)(L_\infty)_*\to\frac{{\rm Lie}(E)(L_\infty)_*}{{\rm Lie}(E)(\sO_L)_*}$ and $\overline\exp_E:{\rm Lie}(E)(L_\infty)_*\to \frac{E(L_\infty)_*}{E(\sO_L)_*}$ the map induced by $\exp_E:{\rm Lie}(E)(L_\infty)\to E(L_\infty)$.
By the commutative diagram
\[\xymatrix{R\Gamma(X,\,\sE_*)[t]\ar[r]\ar@<0.5ex>[d]_{t-A_0~~}\ar@<-0.5ex>[d]
^{~~t-\sum_{s=0}^rA_s\tau^s}&(z^{-d}\widetilde\sO_\infty^n)_*[t]\ar[rr]\ar@<0.5ex>[d]_{t-A_0~~}\ar@<-0.5ex>[d]^{~~t-\sum_{s=0}^rA_s\tau^s}
&&\frac{(L_\infty^n)_*}{(\sO_L^n)_*}[t]
\ar@<0.5ex>[d]_{t-A_0~}\ar@<-0.5ex>[d]^{~t-\sum_{s=0}^rA_s\tau^s}\\
R\Gamma(X,\,\sE_*(e\infty))[t]\ar[r]^\iota&(z^{e-d}\widetilde\sO_\infty^n)_*[t]
\ar[rr]&&\frac{(L_\infty^n)_*}{(\sO_L^n)_*}[t],}\]
we get two quasi-isomorphisms
\begin{eqnarray}\label{f3}
&&R\Gamma\Big(X,\,\sE_*[t]\x{t-A_0}\sE_*(e\infty)[t]\Big)[1]
\simeq\Big(D\x{\bar pp\eta}\frac{{\rm Lie}(E)(L_\infty)_*}{{\rm Lie}(E)(\sO_L)_*}\Big);\\\label{f4}
&&R\Gamma\Big(X,\,\sE_*[t]\x{t-\sum_{s=0}^rA_s\tau^s}\sE_*(e\infty)[t]\Big)[1]
\simeq\Big(D\x{\overline{\exp}_Ep\eta}\frac{E(L_\infty)_*}{E(\sO_L)_*}\Big).
\end{eqnarray}
Since $X$ is projective, then (\ref{f3}) and (\ref{f4}) show that
\begin{eqnarray}\label{li}\Big(D\x{\bar pp\eta}\frac{{\rm Lie}(E)(L_\infty)_*}{{\rm Lie}(E)(\sO_L)_*}\Big)\hbox{ and }\Big(D\x{\overline{\exp}_Ep\eta}\frac{E(L_\infty)_*}{E(\sO_L)_*}\Big)\in D^{\rm perf}(k[t]_*).\end{eqnarray}
By (\ref{f3}) and (\ref{f4}), the composite
\begin{eqnarray*}
&&{\rm det}_{k[t]_*}\Big(R\Gamma\Big(X,\,\sE_*[t]\x{t-A_0}\sE_*(e\infty)[t]\Big)\Big)^{-1}\\
&\simeq&{\rm det}_{k[t]_*}\Big(R\Gamma\Big(X,\,\sE_*[t]\Big)\Big)^{-1}\bigotimes{\rm det}_{k[t]_*}\Big(R\Gamma\Big(X,\,\sE_*(e\infty)[t]\Big)\Big)\\
&\simeq&{\rm det}_{k[t]_*}\Big(R\Gamma\Big(X,\,\sE_*[t]\x{t-\sum_{s=0}^rA_s\tau^s}\sE_*(e\infty)[t]\Big)\Big)^{-1}
\end{eqnarray*}
defines an isomorphism
\begin{eqnarray*}\label{f7}
\gamma:{\rm det}_{k[t]_*}\Big(D\x{\bar pp\eta}\frac{{\rm Lie}(E)(L_\infty)_*}{{\rm Lie}(E)(\sO_L)_*}\Big)
\simeq{\rm det}_{k[t]_*}\Big(D\x{\overline{\exp}_Ep\eta}\frac{E(L_\infty)_*}{E(\sO_L)_*}\Big).
\end{eqnarray*}
By (\ref{f1}), (\ref{f3}), (\ref{f2}) and (\ref{f4}), we get two quasi-isomorphisms
\begin{eqnarray*}\label{f5}
&&g_1:\Big(D\x{\bar pp\eta}\frac{{\rm Lie}(E)(L_\infty)_*}{{\rm Lie}(E)(\sO_L)_*}\Big)
\otimes_{k[t]_*}k((t^{-1}))_*\simeq C;\\\label{f6}
&&g_2:\Big(D\x{\overline{\exp}_Ep\eta}\frac{E(L_\infty)_*}{E(\sO_L)_*}\Big)
\otimes_{k[t]_*}k((t^{-1}))_*\simeq C.
\end{eqnarray*}
By Lemma \ref{c1}, we have
\begin{eqnarray}\label{c2}
L(E,\;*)=[\det(g_2)\circ(\gamma\otimes {\rm id}_{k((t^{-1}))_*})\circ\det(g_1)^{-1}:1]\in k((t^{-1}))_*^\times.
\end{eqnarray}

Let $\sK$ be the kernel of the surjective map $D\x{p\eta}{\rm Lie}(E)(L_\infty)_*$.
%Consider
%\begin{eqnarray}
%&&\label{mu}\Big({\rm Lie}(E)(L_\infty)_*\x{\bar p}\frac{{\rm Lie}(E)(L_\infty)_*}{{\rm Lie}(E)(\sO_L)_*}\Big);\\\label{luo}
%&&\Big({\rm Lie}(E)(L_\infty)_*\x{\overline{\exp}_E}\frac{E(L_\infty)_*}{E(\sO_L)_*}\Big).
%\end{eqnarray}
We get two distinguished triangles
\begin{eqnarray}\label{f8}
&&\sK\to\Big(D\x{\bar pp\eta}\frac{{\rm Lie}(E)(L_\infty)_*}{{\rm Lie}(E)(\sO_L)_*}\Big)\to \Big({\rm Lie}(E)(L_\infty)_*\x{\bar p}\frac{{\rm Lie}(E)(L_\infty)_*}{{\rm Lie}(E)(\sO_L)_*}\Big)\\\label{3174}
&&\sK\to\Big(D\x{\overline{\exp}_Ep\eta}\frac{E(L_\infty)_*}{E(\sO_L)_*}\Big)\to\Big({\rm Lie}(E)(L_\infty)_*\x{\overline{\exp}_E}\frac{E(L_\infty)_*}{E(\sO_L)_*}\Big).
\end{eqnarray}
in $D(k[t]_*)$. Let $L$ be the kernel of the surjective map $C\x{p}{\rm Lie}(E)(L_\infty)_*$.
We get a short exact sequence
\begin{eqnarray}\label{f10}
0\to L\to C\x{p}{\rm Lie}(E)(L_\infty)_*\to0
\end{eqnarray}
of $k((t^{-1}))_*$-modules.
Then $\eta:D\to C$ induces a $k[t]_*$-linear map $\sK\to L$ and a $k((t^{-1}))_*$-linear map $\sK\otimes_{k[t]_*}k((t^{-1}))_*\x{f}L$.
The identity map of ${\rm Lie}(E)(L_\infty)_*$ as $k[t]_*$-modules induces a $k((t^{-1}))_*$-linear map
$${\rm Lie}(E)(L_\infty)_*\otimes_{k[t]_*}k((t^{-1}))_*\to{\rm Lie}(E)(L_\infty)_*$$
and two morphism of complexes
\begin{eqnarray*}
&&h_1:\Big({\rm Lie}(E)(L_\infty)_*\x{\bar p}\frac{{\rm Lie}(E)(L_\infty)_*}{{\rm Lie}(E)(\sO_L)_*}\Big)\otimes_{k[t]_*}k((t^{-1}))_*\to{\rm Lie}(E)(L_\infty)_*;\\
&&h_2:\Big({\rm Lie}(E)(L_\infty)_*\x{\overline{\exp}_E}\frac{E(L_\infty)_*}{E(\sO_L)_*}\Big)\otimes_{k[t]_*}k((t^{-1}))_*\to{\rm Lie}(E)(L_\infty)_*.
\end{eqnarray*}
By \cite{F}, ${\rm Lie}(E)(\sO_L)$ is a lattice of the finite dimensional $k((t^{-1}))$-vector space ${\rm Lie}(E)(L_\infty)$. By the flatness of $k[t]\to k((t^{-1}))$, we have
\begin{eqnarray*}
&&{\rm Hom}_{k[G]}(V,\,{\rm Lie}(E)(\sO_L))\otimes_{F[t]}F((t^{-1}))\\&=&{\rm Hom}_{k[G]}(V,\,{\rm Lie}(E)(\sO_L))\otimes_{k[t]}k((t^{-1}))\\&=&{\rm Hom}_{k[G]}\Big(V,\,{\rm Lie}(E)(\sO_L)\otimes_{k[t]}k((t^{-1}))\Big)\\
&=&{\rm Hom}_{k[G]}(V,\,{\rm Lie}(E)(L_\infty)).
\end{eqnarray*}
This means that ${\rm Hom}_{k[G]}(V,\,{\rm Lie}(E)(\sO_L))$ is a lattice of the finite dimensional $F((t^{-1}))$-vector space ${\rm Hom}_{k[G]}(V,\,{\rm Lie}(E)(L_\infty)).$
So ${\rm Lie}(E)(\sO_L)_*$ is a lattice of the $k((t^{-1}))_*$-module ${\rm Lie}(E)(L_\infty)_*$ for $*=\rho$ or $G$ in the sense of Definition \ref{16}.
Then $h_1$ is a quasi-isomorphism by the fact
$$\Big({\rm Lie}(E)(L_\infty)_*\x{\bar p}\frac{{\rm Lie}(E)(L_\infty)_*}{{\rm Lie}(E)(\sO_L)_*}\Big)\simeq{\rm Lie}(E)(\sO_L)_*\in D^{\rm perf}(k[t]_*).$$
By (\ref{li}) and (\ref{f8}), $\sK\in D^{\rm perf}(k[t]_*)$. By (\ref{li}) and (\ref{3174}),
$$\Big({\rm Lie}(E)(L_\infty)_*\x{\overline{\exp}_E}\frac{E(L_\infty)_*}{E(\sO_L)_*}\Big)\in D^{\rm perf}(k[t]_*).$$
Consider the natural morphism of triangles
\begin{eqnarray}\label{3203}
&&(\ref{f8})\otimes_{k[t]_*}k((t^{-1}))_*\x{(f,g_1,h_1)}(\ref{f10});\\\label{3204}
&&(\ref{3174})\otimes_{k[t]_*}k((t^{-1}))_*\x{(f,g_2,h_2)}(\ref{f10}).
\end{eqnarray}
Since $g_1$ and $h_1$ are quasi-isomorphisms, so is $f$. Since $g_2$ is a quasi-isomorphism, so is $h_2$.
Distinguished triangles (\ref{f8}) and (\ref{3174}) and $\gamma$ define an isomorphism
\begin{eqnarray*}
\gamma_1:{\rm det}_{k[t]_*}\Big({\rm Lie}(E)(L_\infty)_*\x{\bar p}\frac{{\rm Lie}(E)(L_\infty)_*}{{\rm Lie}(E)(\sO_L)_*}\Big)\simeq{\rm det}_{k[t]_*}\Big({\rm Lie}(E)(L_\infty)_*\x{\overline{\exp}_E}\frac{E(L_\infty)_*}{E(\sO_L)_*}\Big).
\end{eqnarray*}
Applying \cite[ Lemma 3.1 ]{F} to (\ref{3203}) and (\ref{3204}), we have
\begin{eqnarray}\label{xiao}
[\det(g_2)(\gamma\otimes{\rm id}_{k((t^{-1}))_*})\det(g_1)^{-1}:1]=[\det(h_2)(\gamma_1\otimes{\rm id}_{k((t^{-1}))_*})\det(h_1)^{-1}:1].
\end{eqnarray}
We have
\begin{eqnarray*}\label{si}
&&H^0\Big({\rm Lie}(E)(L_\infty)_*\x{\overline{\exp}_E}\frac{E(L_\infty)_*}{E(\sO_L)_*}\Big)=\exp_E^{-1}(E(\sO_L))_*;\\\nonumber
&&H^1\Big({\rm Lie}(E)(L_\infty)_*\x{\overline{\exp}_E}\frac{E(L_\infty)_*}{E(\sO_L)_*}\Big)=H(E,*);\\\nonumber
&&H^s\Big({\rm Lie}(E)(L_\infty)_*\x{\overline{\exp}_E}\frac{E(L_\infty)_*}{E(\sO_L)_*}\Big)=0,\hbox{ for } s\neq 0,\;1.
\end{eqnarray*}
Since $h_2$ is a quasi-isomorphism and $\Big({\rm Lie}(E)(L_\infty)_*\x{\overline{\exp}_E}\frac{E(L_\infty)_*}{E(\sO_L)_*}\Big)$ is perfect, then $\exp_E^{-1}(E(\sO_L))_*$ is a lattice of ${\rm Lie}(E)(L_\infty)_*$ and $H(E,\,*)$ is a finite $k[t]_*$-module. Thus we get two lattices ${\rm Lie}(E)(\sO_L)_*$ and $\exp_E^{-1}(E(\sO_L))_*$ of the $k((t^{-1}))_*$-module ${\rm Lie}(E)(L_\infty)_*$.

Fix an isomorphism $\gamma_2:\det_{k[t]_*}\Big({\rm Lie}(E)(\sO_L)_*\Big)\simeq\det_{k[t]_*}\Big(\exp_E^{-1}(E(\sO_L))_*\Big)$ and $\gamma_3:k[t]_*\simeq\det_{k[t]_*}\Big( H(E,\,*)\Big)$.
By \cite[Lemma 3.1]{F} and an equivariant version of \cite [Lemma 3.3 ]{F}, we have
\begin{eqnarray*}
&&[\det(h_2)(\gamma_1\otimes{\rm id}_{k((t^{-1}))_*})\det(h_1)^{-1}:1]\\\nonumber
&=&\frac{[\det\big(H^0(h_2)\big)(\gamma_2\otimes{\rm id}_{k((t^{-1}))_*})\det\big(H^0(h_1)\big)^{-1}:1]}{[\det\big(H^1(h_2)\big)(\gamma_3\otimes{\rm id}_{k((t^{-1}))_*})\det\big(H^1(h_1)\big)^{-1}:1]}\\\nonumber
&=&[{\rm Lie}(E)(\sO_L)_*:\exp_E^{-1}(E(\sO_L))_*]_{k((t^{-1}))_*}\cdot|H(E,\,*)|_{k[t]_*}
\in k((t^{-1}))_*^\times/k_*^\times.
\end{eqnarray*}
Since $L(E,\rho)$, $[{\rm Lie}(E)(\sO_L)_*:\exp_E^{-1}(E(\sO_L))_*]_{k((t^{-1}))_*}$ and $|H(E,\,*)|_{k[t]_*}$ are monic, then by (\ref{c2}) and (\ref{xiao}), we have
$$L(E,\,*)=[{\rm Lie}(E)(\sO_L)_*:\exp_E^{-1}(E(\sO_L))_*]_{k((t^{-1}))_*}\cdot|H(E,\,*)|_{k[t]_*}\in k((t^{-1}))_*^\times.$$
Take $*=\rho$ or $*=G$ if $G$ is an abelian group of order prime to $p$, we have
\begin{eqnarray*}
&&L(E,\,\rho)=[{\rm Hom}_{k[G]}(V,\,{\rm Lie}(E)(\sO_L)):{\rm Hom}_{k[G]}(V,\,\exp_E^{-1}(E(\sO_L)))]_{F((t^{-1}))}\cdot|H(E,\,\rho)|_{F[t]};\\
&&L(E,\,G)=[{\rm Lie}(E)(\sO_L):\exp_E^{-1}(E(\sO_L))]_{k((t^{-1}))[G]}\cdot|H(E,\,G)|_{k[t][G]}.
\end{eqnarray*}
This completes the proofs of Theorem \ref{391} and Theorem \ref{18} (1).
\begin{lem}\label{3101}
For any $k$-module $M$ and any finite $F$-module $N$, the trace map ${\rm tr}_{F/k}:F\to k$ induces ${\rm tr}_{F/k}\otimes{\rm id}_M:F\otimes_kM\to k\otimes_kM=M$ and a $F$-linear isomorphism
$${\rm Hom}_F(N,\,F\otimes_kM)\simeq{\rm Hom}_k(N,\,M);\;f\mapsto {\rm tr}_{F/k}\otimes{\rm id}_{M}\circ f.$$
For any $k[G]$-module $M$ and any finite $F[G]$-module $N$, the trace map ${\rm tr}_{F/k}$ induces an isomorphism
$${\rm Hom}_{F[G]}(N,\,F\otimes_kM)\simeq{\rm Hom}_{k[G]}(N,\,M).$$
\end{lem}
\begin{proof}
Let $0\neq f\in{\rm Hom}_F(N,\,F)$. Then $f$ is surjective. By the surjectivity of ${\rm tr}_{F/k}:F\to k$, ${\rm tr}_{F/k}\circ f:N\to k$ is surjective. This proves
the injectivity of ${\rm Hom}_F(N,\,F)\x{{\rm tr}_{F/k}}{\rm Hom}_k(N,\,k)$. Both sides are $F$-vector spaces of the same dimension, then ${\rm tr}_{F/k}:{\rm Hom}_F(N,\,F)\to{\rm Hom}_k(N,\,k)$ is an isomorphism of $F$-vector spaces. The lemma follows from the canonical isomorphisms
$${\rm Hom}_F(N,\,F)\otimes_F(F\otimes_kM)\simeq{\rm Hom}_F(N,\,F\otimes_kM)\hbox{ and }{\rm Hom}_k(N,\,k)\otimes_kM\simeq{\rm Hom}_k(N,\,M).$$

\end{proof}
\begin{lem}\label{3102}
Let $N$ be a finite $F[G]$-module and let $\phi$ be a $k[t][G]$-linear endomorphism of a finitely generated free $k[t][G]$-module $M$. Let $\phi^*$ be the $F[t]$-linear endomorphism of the $F[t]$-module ${\rm Hom}_{k[G]}(N,\,M)$ induced by $\phi$. Suppose $G$ is an abelian group of order prime to $p$.
We have
$${\rm det}_{F[t]}\Big(\phi^*,{\rm Hom}_{k[G]}(N,\,M)\Big)={\rm det}_{F[t]}\Big({\rm det}_{k[t][G]}(\phi,\,M),\;N[t]\Big).$$
\end{lem}
\begin{proof}
By Lemma \ref{3101}, we have
\begin{eqnarray*}
&&{\rm det}_{F[t]}\Big(\phi^*,{\rm Hom}_{k[G]}(N,\,M)\Big)={\rm det}_{F[t]}\Big(({\rm id}_F\otimes\phi)^*,{\rm Hom}_{F[G]}(N,\,F\otimes_kM)\Big);\\
&&{\rm det}_{F[t]}\Big({\rm det}_{k[t][G]}(\phi,\,M),\;N[t]\Big)={\rm det}_{F[t]}\Big({\rm det}_{F[t][G]}({\rm id}_F\otimes\phi,\,F\otimes_kM),\;N[t]\Big).
\end{eqnarray*}
So we may assume that $F=k$ which is algebraically closed. For any character $\chi:G\to k^\times$, set $e_\chi=\sum_{g\in G}\frac{1}{|G|}\chi(g)^{-1}g\in k[G]$. Then $e_\chi$ is an idempotent element of $k[G]$ and $e_\chi k[G]\simeq k$ and
$k[G]=\prod_{\chi:G\to k^\times}e_\chi k[G]$.
Then
\begin{eqnarray*}
&&{\rm det}_{k[t]}\Big(\phi^*,{\rm Hom}_{k[G]}(N,\,M)\Big)=\prod_{\chi:G\to k^\times}{\rm det}_{k[t]}\Big((e_\chi\phi)^*,{\rm Hom}_{e_\chi k[G]}(e_\chi N,\,e_\chi M)\Big);\\
&&{\rm det}_{k[t]}\Big({\rm det}_{k[G][t]}(\phi,\,M),\;N[t]\Big)=\prod_{\chi:G\to k^\times}{\rm det}_{k[t]}\Big({\rm det}_{e_\chi k[G][t]}(e_\chi\phi,\,e_\chi M),\;e_\chi N[t]\Big)
\end{eqnarray*}
Then we may assume that $G$ is the trivial group. In this case,
$${\rm det}_{k[t]}\Big(\phi^*,{\rm Hom}_{k}(N,\,M)\Big)={\rm det}_{k[t]}(\phi,\,M)^{\dim_kN}={\rm det}_{k[t]}\Big({\rm det}_{k[t]}(\phi,\,M),\;N[t]\Big).$$
%Since $G$ is abelian of order prime to $p$, then $k[G]$
%Any finite $k[G]$-module is semisimple and any finite simple $k[G]$-module is of dimension one. Then $N$ is a one dimensional $k$-representation of $G$ of character $\chi:G\to k^\times$. L. Then $e_\chi^2=e_\chi$ and $ge_\chi=\chi(g)e_\chi$ for any $g\in G$. Choose a $k[t][G]$-basis of $\{m_1,\ldots,m_t\}$ of $M$ and let $(a_{ij})_{ij}$ be the matrix of $\phi$ on $M$ with respect to this basis. Then $\{e_\chi m_1,\ldots,e_\chi m_t\}$ is a $k[t]$-basis of ${\rm Hom}_{k[G]}(N,\,M)=e_\chi M$ and $(\chi(a_{ij}))_{ij}$ is the matrix of $\phi_*$ on $e_\chi M$.
%Then $${\rm det}_{k[t]}(\phi_*,\,e_\chi M)=\det(\chi(a_{ij}))=\chi(\det(a_{ij}))={\rm det}_{k[t]}(\det(a_{ij}),\,N)={\rm det}_{k[t]}\Big({\rm det}_{k[t][G]}(\phi,\,M),\,N\Big).$$
\end{proof}
Applying Lemma \ref{3102} to $N=V$, $M=\sO_L^n/\mathfrak p\sO_L^n[t]$ and $\phi=t-A_0$ or $\phi=t-\sum_{s=0}^rA_s\tau^s$, we have
\begin{eqnarray*}
L(E,\,\rho)&=&\prod_{\mathfrak p\in{\rm Max}(\sO_K)}\frac{{\rm det}_{F[t]}\Big(t-A_0,\;{\rm Hom}_{k[G]}(V,\,\sO_L^n/\mathfrak p\sO_L^n[t])\Big)}{{\rm det}_{F[t]}
\Big(t-\sum_{s=0}^rA_s\tau^s,\;{\rm Hom}_{k[G]}(V,\,\sO_L^n/\mathfrak p\sO_L^n[t])\Big)}\\
&=&\prod_{\mathfrak p\in{\rm Max}(\sO_K)}\frac{{\rm det}_{F[t]}\Big({\rm det}_{k[t][G]}\Big(t-A_0,\,\sO_L^n/\mathfrak p\sO_L^n[t]\Big),\,V[t]\Big)}
{{\rm det}_{F[t]}\Big({\rm det}_{k[t][G]}\Big(t-\sum_{s=0}^rA_s\tau^s,\,\sO_L^n/\mathfrak p\sO_L^n[t]\Big),\,V[t]\Big)}\\
&=&{\rm det}_{F((t^{-1}))}\Big(\prod_{\mathfrak p\in{\rm Max}(\sO_K)}\frac{{\rm det}_{k[t][G]}\Big(t-A_0,\,\sO_L^n/\mathfrak p\sO_L^n[t]\Big)}
{{\rm det}_{k[t][G]}\Big(t-\sum_{s=0}^rA_s\tau^s,\,\sO_L^n/\mathfrak p\sO_L^n[t]\Big)},\,V((t^{-1}))\Big)\\
&=&{\rm det}_{F((t^{-1}))}\Big(L(E,\,G),\,V((t^{-1}))\Big)\in F((t^{-1})).
\end{eqnarray*}
This proves Theorem \ref{18} (2).
\section{Special values of Artin $L$-functions of Galois representations}
Recall that $L/K$ is a finite Galois extension of Galois group $G$. Let $\mathfrak p$ be a maximal ideal of $\sO_K$ and $\mathfrak P$ a maximal ideal of $\sO_L$ above $\mathfrak p$. Let $e$ be the ramification index of $L/K$ at $\mathfrak p$. Let ${\rm frob}_\mathfrak P$ be the inverse image of the Frobenius element ${\rm Frob}_\mathfrak P\in G_\mathfrak P/I_\mathfrak P $ in $G_\mathfrak P$. Then $|{\rm frob}_\mathfrak P|=e$.
\begin{lem}\label{3131}
If $L/K$ is abelian and tamely ramified at $\mathfrak p$, then $e\big||\kappa_\mathfrak p^\times|$.
\end{lem}
\begin{proof}
Let $K_\mathfrak p$ and $L_\mathfrak P$ be the complete fields of $K$ and $L$ at $\mathfrak p$ and $\mathfrak P$, respectively. Then $G_\mathfrak P\simeq G(L_\mathfrak P/K_\mathfrak p)$ and $L_\mathfrak P/K_\mathfrak p$ is tamely ramified at $\mathfrak p$. Let $M$ be the maximal unramified extension of $K_\mathfrak p$ in $L_\mathfrak P$. So $L_\mathfrak P/M$ is totally and tamely ramified. So there exists a uniformizer $\xi$ of $M$ such that $L_\mathfrak P=M(\sqrt[e]\xi)$. Since $L_\mathfrak P/M$ is Galois, one can choose a primitive $e$-th root of unity $\zeta$ in $M$. Let $\sigma\in G(L_\mathfrak P/M)$ such that $\sigma(\sqrt[e]\xi)=\zeta\sqrt[e]\xi.$ Choose $\delta\in G(L_\mathfrak P/K_\mathfrak p)$ such that $\delta|_M$ is the Frobenius element of the unramified extension $M/K_\mathfrak p$. Then $\delta(\sqrt[e]\xi)=\alpha\sqrt[e]\xi$ for some $\alpha\in\sO_{L_\mathfrak P}^\times$.
One has
$$\sigma(\delta(\sqrt[e]\xi))=\sigma(\alpha\sqrt[e]\xi)=\sigma(\alpha)\zeta\sqrt[e]\xi\hbox{ and }\delta(\sigma(\sqrt[e]\xi))=\delta(\zeta\sqrt[e]\xi)=\delta(\zeta)\alpha\sqrt[e]\xi.$$
Since $G$ is abelian, we have $\sigma(\alpha)\zeta=\delta(\zeta)\alpha$. Since $L_\mathfrak P/M$ is totally ramified, then the residue fields of $M$ equal to $\kappa_\mathfrak P$.
Let $\bar\delta\in G(\kappa_\mathfrak P/\kappa_\mathfrak p)$ defined by $\delta$ and $\bar\zeta$ the image of $\zeta$ in $\kappa_\mathfrak P$. Then $\bar\delta(\bar\zeta)=\bar\zeta$. Since $\bar\delta$ generates $G(\kappa_\mathfrak P/\kappa_\mathfrak p)$, then $\bar\zeta\in\kappa_\mathfrak p.$ Since $(e,\,|\kappa_\mathfrak p|)=1$, then $\bar\zeta$ is a primitive $e$-th root of unity in $\kappa_\mathfrak p$ and hence $e\big||\kappa_\mathfrak p^\times|$.
\end{proof}

\begin{lem}\label{3132} Suppose $L/K$ is abelian and tamely ramified at $\mathfrak p$.
Then $$\sum_{\sigma\in{\rm frob}_\mathfrak P}\sigma(x)\equiv ex^{|\kappa_\mathfrak p|}\;({\rm mod}\; \mathfrak p\sO_L)\hbox{ for any }x\in\sO_L.$$
\end{lem}
\begin{proof}
Let $M=L^{I_\mathfrak P}$ and let $\mathfrak q_1,\ldots,\mathfrak q_r$ be the various prime ideals of $\sO_M$ above $\mathfrak p$. Then $\mathfrak p\sO_M=\mathfrak q_1\cdots\mathfrak q_r$. There exists a unique prime ideal $\mathfrak P_i$ of $\sO_L$ above $\mathfrak q_i$ for any $1\leq i\leq r$. Then $\mathfrak q_i\sO_L=\mathfrak P_i^e$ and $\sO_M/\mathfrak q_i\simeq \sO_L/\mathfrak P_i$. For any $x\in\sO_L$, there exists $y\in\sO_M$ and $z\in\mathfrak P_i$ such that $x=y+z$. By the definition of ${\rm frob}_\mathfrak P$, we have
$$\sum_{\sigma\in{\rm frob}_\mathfrak P}(\sigma(y)-y^{|\kappa_\mathfrak p|})\in\mathfrak q_i.$$
By Lemma \ref{3131}, $e<|\kappa_\mathfrak p|$ and then
$$\sum_{\sigma\in{\rm frob}_\mathfrak P}\sigma(z)\in L^{I_\mathfrak P}\cap\mathfrak P_i=M\cap\mathfrak P_i=\mathfrak q_i\hbox{ and }z^{|\kappa_\mathfrak p|}\in\mathfrak P_i^e=\mathfrak q_i\sO_L.$$
So
$$\sum_{\sigma\in{\rm frob}_\mathfrak P}(\sigma(x)-x^{|\kappa_\mathfrak p|})=\sum_{\sigma\in{\rm frob}_\mathfrak P}(\sigma(y)-y^{|\kappa_\mathfrak p|})+\sum_{\sigma\in{\rm frob}_\mathfrak P}(\sigma(z)-z^{|\kappa_\mathfrak p|})\in\bigcap_{i=1}^r\mathfrak q_i\sO_L=\mathfrak p\sO_L.$$
\end{proof}

\begin{lem}\label{3133}
Suppose $L/K$ is abelian and tamely ramified at $\mathfrak p$. Then
\begin{eqnarray*}
&&|C^{\otimes n}(\sO_L/\mathfrak p\sO_L)|_{k[t][G]}=|\kappa_\mathfrak p|^n_{k[t]}-\frac{1}{e}\sum_{\sigma\in{\rm frob}_\mathfrak P}\sigma\in k[t][G];\\
&&|{\rm Lie}(C^{\otimes n})(\sO_L/\mathfrak p\sO_L)|_{k[t][G]}=|\kappa_\mathfrak p|^n_{k[t]}\in k[t].
\end{eqnarray*}
\end{lem}
\begin{proof}
%Let $f$ be the monic generator of $\mathfrak p\cap k[t]$. Then $|\kappa_\mathfrak p|_{k[t]}=f^{[\kappa_\mathfrak p:\kappa_{(f)}]}$.
%Write $f=\prod_{i\in\Z/r\Z}(t-t_i)$ with elements $t_i$ in the algebraic closure $\bar k$ of $k$ such that $t_{i+1}=t_i^{q}$ for $i\in\Z/r\Z$.
%By Definition \ref{3134}, $C^\otimes n$ is defined by $A_0=zI_n+N$ and $A_1$.
%Let
%$$M_i=\{x\in\sO_L^n/\mathfrak p\sO_L^n\otimes_k\bar k|\;tx=t_ix\}\hbox{ for }\;1\leq i\leq r.$$
%Then $M_i$ is free $\bar k$-module and $A_0(M_i)\subset M_i$ and $A_1\tau(M_i)\subset M_{i-1}$ for any $i$.
For any $1\leq i\leq n$ and $\delta\in G(\kappa_\mathfrak p/k)$, let $M_{i,\delta}$ be the $k$-vector space $\sO_L/\mathfrak p\sO_L$ equipped with
the $\kappa_\mathfrak p$-module structure given by the composite
$\kappa_\mathfrak p\x{\delta}\kappa_\mathfrak p=\sO_K/\mathfrak p\hookrightarrow\sO_L/\mathfrak p\sO_L.$
Then we have a $\kappa_\mathfrak p$-linear isomorphism
$$\sO_L^n/\mathfrak p\sO_L^n\otimes_k\kappa_\mathfrak p\simeq\bigoplus_{1\leq i\leq n,\;\delta\in G(\kappa_\mathfrak p/k)}M_{i,\delta};\;(x_1,\ldots,x_n)\otimes y\mapsto\oplus_{i,\delta}\delta(y)x_i,$$
where the $\kappa_\mathfrak p$-module structure on $\sO_L^n/\mathfrak p\sO_L^n\otimes_k\kappa_\mathfrak p$ is given by
$\kappa_\mathfrak p$. By Definition \ref{3151},
$$C^{\otimes n}(t)(x_1,\ldots,x_{n-1},x_n)=(zx_1,\ldots,zx_{n-1},zx_n)+(x_2,\ldots,x_n,x_1^q)$$ for any $(x_1,\ldots,x_n)\in\sO_L^n/\mathfrak p\sO_L^n$.
Let $\tau\in G(\kappa_\mathfrak p/k)$ be the $q$-th power map.
Then the operator $z\otimes {\rm id}_{\kappa_\mathfrak p}$ on $\sO_L^n/\mathfrak p\sO_L^n\otimes_k\kappa_\mathfrak p$ corresponds to
$\oplus_{i,\delta}\delta^{-1}(z):\bigoplus_{i,\delta}M_{i,\delta}\to\bigoplus_{i,\delta} M_{i,\delta}$ and that of $(C^{\otimes n}(t)-z)\otimes{\rm id}_{\kappa_\mathfrak p}$ to
$\phi$ such that
$\phi(M_{i,\delta})\subset M_{i-1,\delta}$ and $\phi(M_{1,\delta})\subset M_{n,\delta\tau}$ for any $2\leq i\leq n$ and $\delta\in G(\kappa_\mathfrak p/k)$. Moreover,
$\phi:M_{i,\delta}\to M_{i-1,\delta}$ corresponds to the identity map and $\phi:M_{1,\delta}\to M_{n,\delta\tau}$ to the $q$-th power map
via $M_{i,\delta}=\sO_L/\mathfrak p\sO_L$. Then $\phi^{n[\kappa_\mathfrak p:k]}$ preserves $M_{i,\delta}$ and $\phi^{n[\kappa_\mathfrak p:k]}:M_{i,\delta}\to M_{i,\delta}$ corresponds to the $q^{[\kappa_\mathfrak p:k]}=|\kappa_\mathfrak p|$-th power map via $M_{i,\delta}=\sO_L/\mathfrak p\sO_L$.
By Lemma \ref{3132}, $\phi^{n[\kappa_\mathfrak p:k]}=\frac{1}{e}\sum_{\sigma\in{\rm frob}_\mathfrak P}\sigma$ on $M_{i,\delta}$.
Since $L/K$ is tamely ramified at $\mathfrak p$, then $\sO_L/\mathfrak p\sO_L$ is a
free $\kappa_\mathfrak p[G]$-module of rank one.
By a little modification of \cite[ Lemma 8.1.3 ]{BP}, we have
\begin{eqnarray*}
&&|C^{\otimes n}(\sO_L/\mathfrak p\sO_L)|_{k[t][G]}\\
&=&{\rm det}_{k[t][G]}\Big(t-C^{\otimes n}(t),\;\sO_L^n/\mathfrak p\sO_L^n[t]\Big)\\
&=&{\rm det}_{\kappa_\mathfrak p[t][G]}\Big(t-C^{\otimes n}(t)\otimes{\rm id}_{\kappa_\mathfrak p},\;\sO_L^n/\mathfrak p\sO_L^n\otimes_k\kappa_\mathfrak p[t]\Big)\\
&=&{\rm det}_{\kappa_\mathfrak p[t][G]}\Big(t-\oplus_{i,\,\delta}\delta^{-1}(z)-\phi,\;\bigoplus_{1\leq i\leq n,\,\delta\in G(\kappa_\mathfrak p/k)}M_{i,\delta}[t]\Big)\\
&=&{\rm det}_{\kappa_\mathfrak p[t][G]}\Big(\prod_{1\leq i\leq n,\,\delta\in G(\kappa_\mathfrak p/k)}(t-\delta^{-1}(z))-\phi^{n[\kappa_\mathfrak p:k]},\;M_{1,\tau}[t]\Big)\\
&=&\prod_{1\leq i\leq n,\;\delta\in G(\kappa_\mathfrak p/k)}(t-\delta^{-1}(z))-\frac{1}{e}\sum_{\sigma\in{\rm frob}_\mathfrak P}\sigma\\
&=&|\kappa_\mathfrak p|_{k[t]}^n-\frac{1}{e}\sum_{\sigma\in{\rm frob}_\mathfrak P}\sigma\in k[t][G].
\end{eqnarray*}
Using the same method, we also have
$$|{\rm Lie}(C^{\otimes n})(\sO_L/\mathfrak p\sO_L)|_{k[t][G]}=|\kappa_\mathfrak p|^n_{k[t]}\in k[t].$$
\end{proof}
By normal basis theorem, $\kappa_\mathfrak P$ is a free $\kappa_\mathfrak p[G(\kappa_\mathfrak P/\kappa_\mathfrak p)]$-module of rank one. By the same argument of Lemma \ref{3133}, we have the following corollary for any finite Galois extension $L/K$ not necessary abelian or tamely ramified at $\mathfrak p$.
\begin{cor}\label{3152}
We have
\begin{eqnarray*}
&&|C^{\otimes n}(\kappa_\mathfrak P)|_{k[t][G(\kappa_\mathfrak P/\kappa_\mathfrak p)]}=|\kappa_\mathfrak p|^n_{k[t]}-{\rm Frob}_\mathfrak P\in k[t][G(\kappa_\mathfrak P/\kappa_\mathfrak p)];\\\nonumber
&&|{\rm Lie}(C^{\otimes n})(\kappa_\mathfrak P)|_{k[t][G(\kappa_\mathfrak P/\kappa_\mathfrak p)]}=|\kappa_\mathfrak p|^n_{k[t]}\in k[t].
\end{eqnarray*}
\end{cor}

\begin{lem}\label{3242}
Let $\rho:G\to{\rm GL}_F(V)$ be a finite dimensional $F$-linear representation. We have
\begin{eqnarray}
\label{shuang}\frac{|{\rm Hom}_{k[G]}(V,\,C^{\otimes n}(\sO_L/\mathfrak p\sO_L))|_{F[t]}}{|{\rm Hom}_{k[G]}(V,\,{\rm Lie}(C^{\otimes n})(\sO_L/\mathfrak p\sO_L))|_{F[t]}}
={\rm det}_{F[[t^{-1}]]}\Big(1-\frac{\rho({\rm Frob}_\mathfrak P)}{|\kappa_\mathfrak p|_{k[t]}^n},\,V_{I_\mathfrak P}[[t^{-1}]]\Big);\\\label{3154}
\frac{|V^*\otimes_{k[G]}C^{\otimes n}(\sO_L/\mathfrak p\sO_L)|_{F[t]}}{|V^*\otimes_{k[G]}{\rm Lie}(C^{\otimes n})(\sO_L/\mathfrak p\sO_L)|_{F[t]}}={\rm det}_{F[[t^{-1}]]}\Big(1-\frac{\rho({\rm Frob}_\mathfrak P)}{|\kappa_\mathfrak p|_{k[t]}^n},\,V^{I_\mathfrak P}[[t^{-1}]]\Big).
\end{eqnarray}
\end{lem}
\begin{proof}
Let $\mathfrak P_1,\ldots,\mathfrak P_r$ be the various primes of $\sO_L$ lying above $\mathfrak p$ such that $\mathfrak P_1=\mathfrak P$. Then $\sO_L/\mathfrak p\sO_L=\prod_{i=1}^r\sO_L/\mathfrak P_i^e$. So $\sO_L/\mathfrak p\sO_L={\rm Ind}^G_{G_\mathfrak P}(\sO_L/\mathfrak P^e)=k[G]\otimes_{k[G_\mathfrak P]}\sO_L/\mathfrak P^e$ as $k[G]$-modules. Then
\begin{eqnarray*}
&&{\rm Hom}_{k[G]}(V,\,\sO_L/\mathfrak p\sO_L)={\rm Hom}_{k[G]}(V,\,{\rm Ind}_{G_\mathfrak P}^G(\sO_L/\mathfrak P^e))={\rm Hom}_{k[G_\mathfrak P]}(V,\,\sO_L/\mathfrak P^e);\\
&&V^*\otimes_{k[G]}\sO_L/\mathfrak p\sO_L\simeq V^*\otimes_{k[G]}k[G]\otimes_{k[G_\mathfrak P]}\sO_L/\mathfrak P^e=V^*\otimes_{k[G_\mathfrak P]}\sO_L/\mathfrak P^e.
\end{eqnarray*}
So we may assume $r=1$ and then $G=G_\mathfrak P$. Let $\tau_0={\rm Lie}(C^{\otimes n})(t)$ and $\tau_1=C^{\otimes n}(t)-{\rm Lie}(C^{\otimes n})(t)$.
Since $((\mathfrak P/\mathfrak P^e)^n,\,\tau_0,\tau_1)\in{\rm NilCoh}_{\widetilde\tau}^G(\kappa_\mathfrak p,\,k)$, the short exact sequence $0\to\mathfrak P/\mathfrak P^{e}\to\sO_L/\mathfrak P^{e}\to\kappa_\mathfrak P\to0$ shows
\begin{eqnarray}\label{3176}
({\rm Hom}_{k[G]}(V,\,(\sO_L/\mathfrak P^e)^n),\tau_0,\,\tau_1)\simeq({\rm Hom}_{k[G]}(V,\,\kappa_\mathfrak P^n),\,\tau_0,\,\tau_1)\in\widetilde{\rm Crys}(\kappa_\mathfrak p,\,F).
\end{eqnarray}
%%$$0\to{\rm Hom}_{k[G]}(V,\,\mathfrak P^{i-1}/\mathfrak P^i)\to {\rm Hom}_{k[G]}(V,\,\sO_L/\mathfrak P^i)\to {\rm Hom}_{k[G]}(V,\,\sO_L/\mathfrak P^{i-1})\x{\delta}{\rm Ext}_{k[G]}^1(V,\,\mathfrak P^{i-1}/\mathfrak P^i).$$
%Since the $q$-th power map on $\mathfrak P^{i-1}/\mathfrak P^i$ is zero and hence $C^{\otimes n}(t)={\rm Lie}(C^{\otimes n})(t)$ on ${\rm Hom}_{k[G]}(V,\,\mathfrak %P^{i-1}/\mathfrak P^i)$ and ${\rm im}(\delta).$ Then
%$${\rm det}_{F[t]}\Big(t-C^{\otimes n}(t),\,{\rm Hom}_{k[G]}(V,\,\mathfrak P^{i-1}/\mathfrak P^i)[t]\Big)={\rm det}_{F[t]}\Big(t-{\rm Lie}(C^{\otimes n})(t),\,{\rm %Hom}_{k[G]}(V,\,\mathfrak P^{i-1}/\mathfrak P^i)[t]\Big).$$
%For $2\leq i\leq e$, we have
%\begin{eqnarray*}
%\frac{{\rm det}_{F[t]}\Big(t-C^{\otimes n}(t),\,{\rm Hom}_{k[G]}(V,\,\sO_L/\mathfrak P^i)[t]\Big)}{{\rm det}_{F[t]}\Big(t-{\rm Lie}(C^{\otimes n})(t),\,{\rm %Hom}_{k[G]}(V,\,\sO_L/\mathfrak P^i)[t]\Big)}
%=\frac{{\rm det}_{F[t]}\Big(t-C^{\otimes n}(t),\,{\rm Hom}_{k[G]}(V,\,\sO_L/\mathfrak P^{i-1})[t]\Big)}{{\rm det}_{F[t]}\Big(t-{\rm Lie}(C^{\otimes n})(t),\,{\rm %Hom}_{k[G]}(V,\,\sO_L/\mathfrak P^{i-1})[t]\Big)}.
%\end{eqnarray*}
%Then
Since $I_\mathfrak P$ acts trivially on $\kappa_\mathfrak P$ and $G/I_\mathfrak P\simeq G(\kappa_\mathfrak P/\kappa_\mathfrak p)$, we have
\begin{eqnarray*}
&&{\rm Hom}_{k[G]}(V,\, \kappa_\mathfrak P)\simeq{\rm Hom}_{k[G]}(V_{I_\mathfrak P},\,\kappa_\mathfrak P)={\rm Hom}_{k[G(\kappa_\mathfrak P/\kappa_\mathfrak p)]}(V_{I_\mathfrak P},\,\kappa_\mathfrak P)\\
&&V^*\otimes_{k[G]}\kappa_\mathfrak P=(V^*)_{I_\mathfrak P}\otimes_{k[G]}\kappa_\mathfrak P\simeq (V^{I_\mathfrak P})^*\otimes_{k[G(\kappa_\mathfrak P/\kappa_\mathfrak p)]}\kappa_\mathfrak P.
\end{eqnarray*}

By Lemma \ref{3102} and Corollary \ref{3152}, we have
\begin{eqnarray*}
&&\frac{|{\rm Hom}_{k[G]}(V,\,C^{\otimes n}(\sO_L/\mathfrak p\sO_L))|_{F[t]}}{|{\rm Hom}_{k[G]}(V,\,{\rm Lie}(C^{\otimes n})(\sO_L/\mathfrak p\sO_L))|_{F[t]}}\\
&=&\frac{{\rm det}_{F[t]}\Big(t-C^{\otimes n}(t),\,{\rm Hom}_{k[G]}(V,\,\sO_L^n/\mathfrak p\sO_L^n)[t]\Big)}{{\rm det}_{F[t]}\Big(t-{\rm Lie}(C^{\otimes n})(t),\,{\rm Hom}_{k[G]}(V,\,\sO_L^n/\mathfrak p\sO_L^n)[t]\Big)}\\
&=&\frac{{\rm det}_{F[t]}\Big(t-C^{\otimes n}(t),\,{\rm Hom}_{k[G(\kappa_\mathfrak P/\kappa_\mathfrak p)]}(V_{I_\mathfrak P},\,\kappa_\mathfrak P^n[t])\Big)}{{\rm det}_{F[t]}\Big(t-{\rm Lie}(C^{\otimes n})(t),\,{\rm Hom}_{k[G(\kappa_\mathfrak P/\kappa_\mathfrak p)]}(V_{I_\mathfrak P},\,\kappa_\mathfrak P^n[t])\Big)}\\
&=&\frac{{\rm det}_{F[t]}\Big({\rm det}_{k[G(\kappa_\mathfrak P/\kappa_\mathfrak p)][t]}\Big(t-C^{\otimes n}(t),\;\kappa_\mathfrak P^n[t]\Big),\;V_{I_\mathfrak P}[t]\Big)}{{\rm det}_{F[t]}\Big({\rm det}_{k[G(\kappa_\mathfrak P/\kappa_\mathfrak p)][t]}\Big(t-{\rm Lie}(C^{\otimes n})(t),\;\kappa_\mathfrak P^n[t]\Big),\;V_{I_\mathfrak P}[t]\Big)}\\
&=&\frac{{\rm det}_{F[t]}\Big(|\kappa_\mathfrak p|_{k[t]}^n-\rho({\rm Frob}_\mathfrak P),\;V_{I_\mathfrak P}[t]\Big)}{{\rm det}_{F[t]}\Big(|\kappa_\mathfrak p|_{k[t]}^n,\;V_{I_\mathfrak P}[t]\Big)}\\
&=&{\rm det}_{F[[t^{-1}]]}\Big(1-\frac{\rho({\rm Frob}_\mathfrak P)}{|\kappa_\mathfrak p|^n_{k[t]}},\,V_{I_\mathfrak P}[[t^{-1}]]\Big)\in 1+t^{-1}F[[t^{-1}]],
\end{eqnarray*}
where the second equality holds by applying Lemma \ref{3172} to (\ref{3176}).

Equality (\ref{3154}) is a dual version of (\ref{shuang}).
\end{proof}
Define the trace map
$${\rm tr_\mathfrak P}:V_{I_\mathfrak P}\to V^{I_\mathfrak P},\;{\rm tr}_\mathfrak P(v)=\sum_{g\in I_\mathfrak P}\rho(g)(v)\hbox{ for any }v\in V_{I_\mathfrak P}.$$
Suppose $L/K$ is tamely ramified at $\mathfrak p$. This means $(p,\,e)=1$. For any $v\in V^{I_\mathfrak P}$, we have ${\rm tr}_\mathfrak P(\frac{1}{e}v)=v$.
For any $v\in V$ such that ${\rm tr}_\mathfrak P(v)=0$, we have $v=\sum_{g\in I_\mathfrak P}\frac{1}{e}(v-\rho(g)(v))\in {_{I_\mathfrak P}}V$. This proves that ${\rm tr}_\mathfrak P$ is an isomorphism if $L/K$ is tamely ramified at $\mathfrak p$. Since $L/K$ is Galois, there are only finitely many prime ideals of $\sO_K$ where $L/K$ is ramified.
By Lemma \ref{3242}, we prove Lemma \ref{3164} and Lemma \ref{3241}.
%$$L(n,\,\rho)=\prod_{\mathfrak p\in{\rm Max}(\sO_K)}\frac{|V^*\otimes_{k[G]}{\rm Lie}(C^{\otimes n})(\sO_L/\mathfrak p\sO_L)|_{F[t]}}{|V^*\otimes_{k[G]}C^{\otimes n}(\sO_L/\mathfrak p\sO_L)|_{F[t]}}\in F((t^{-1})).$$

\begin{lem}\label{zhi}
Define the trace map
$${\rm tr}_{\rho}:V^*\otimes_{k[G]}\sO_L^n\to{\rm Hom}_{k[G]}(V,\,\sO_L^n)$$
by ${\rm tr}_\rho(f\otimes x)(v)=\sum_{g\in G}{\rm tr}_{F/k}(f(g^{-1}v))gx$ for any $f\in V^*,\,x\in\sO_L^n$ and $v\in V$.

(1) The trace map ${\rm tr}_\rho$ is well defined and $F\otimes_k\sO_K$-linear.

(2) Let $S$ be the finite set of prime ideals of $\sO_K$ where $L/K$ is wildly ramified. Let $M={\rm ker}({\rm tr}_\rho)$ and $N={\rm cok}({\rm tr}_\rho)$. Then $M$ and $N$ are two finite $\sO_K$-modules supported on $S$.

(3) If $V$ is a projective $F[G]$-module, then $\alpha$ is an isomorphism.
\end{lem}
\begin{proof}
(1) holds by direct calculations. Let $\sO_{K,\,\mathfrak p}$ be the complete local ring of $\sO_K$ at $\mathfrak p$ and $\sO_{L,\,\mathfrak p}=\sO_L\otimes_{\sO_K}\sO_{K,\,\mathfrak p}$. Suppose $L/K$ is tamely ramified at $\mathfrak p$. Then $\sO_{L,\,\mathfrak p}$ is a free $\sO_{K,\,\mathfrak p}[G]$-module of rank one. To prove (2), it suffices to prove the bijectivity of
$${\rm tr}_{\rho}:V^*\otimes_{k[G]}\sO_{L,\,\mathfrak p}^n\to{\rm Hom}_{k[G]}(V,\,\sO_{L,\,\mathfrak p}^n).$$
Since $\sO_{L,\,\mathfrak p}$ is a free $k[G]$-module, then we only need to show the bijectivity of
$${\rm tr}_{\rho}:V^*\otimes_{k[G]}k[G]\to{\rm Hom}_{k[G]}(V,\,k[G]).$$
The proof is direct. If $V$ is a projective $F[G]$-module, we may assume that $V=F[G]$ is the regular representation ${\rm reg}$ of $G$.
Then (3) holds by the isomorphism
$${\rm tr}_{\rm reg}:{\rm Hom}_{F}(F[G],\,F)\otimes_{k[G]}\sO_L^n\to{\rm Hom}_{k[G]}(F[G],\,\sO_L^n).$$
\end{proof}
\begin{thm}
We have
\begin{eqnarray*}
L(n,\,\rho)&=&[{\rm Hom}_{k[G]}(V,\,{\rm Lie}(C^{\otimes n})(\sO_L)):{\rm Hom}_{k[G]}(V,\,\exp_{C^{\otimes n}}^{-1}(C^{\otimes n}(\sO_L)))]_{F((t^{-1}))}\\
&&\cdot|H(C^{\otimes n},\,\rho)|_{F[t]}\cdot\frac{|{\rm Lie}(C^{\otimes n})(M)|_{F[t]}}{|C^{\otimes n}(M)|_{F[t]}}\cdot\frac{|C^{\otimes n}(N)|_{F[t]}}{|{\rm Lie}(C^{\otimes n})(N)|_{F[t]}}.
\end{eqnarray*}
This means that $L(C^{\otimes n},\,\rho)=\lambda\cdot L(n,\,\rho)$ for some $\lambda\in F(t)^\times$. If $L/K$ is tamely ramified at any prime ideal of $\sO_K$ or $V$ is a projective $F[G]$-module, then $L(C^{\otimes n},\,\rho)=L(n,\,\rho)$.
\end{thm}
\begin{proof}
By Lemma \ref{zhi} and applying Lemma \ref{3172} to the exact sequence
$$0\to M\to V^*\otimes_{k[G]}\sO_L^n\x{{\rm tr}_\rho}{\rm Hom}_{k[G]}(V,\,\sO_L^n)\to  N\to 0,$$
we have
\begin{eqnarray*}
L(n,\,\rho)=L(C^{\otimes n},\,\rho)\cdot\frac{|{\rm Lie}(C^{\otimes n})(M)|_{F[t]}}{|C^{\otimes n}(M)|_{F[t]}}\cdot\frac{|C^{\otimes n}(N)|_{F[t]}}{|{\rm Lie}(C^{\otimes n})(N)|_{F[t]}}.
\end{eqnarray*}
The theorem holds by applying Theorem \ref{391} to $C^{\otimes n}$.
\end{proof}

%Thus
%\begin{eqnarray*}
%&&R\Gamma\Big(X,\,\sO_X(-d\infty)^n[t]\x{t-A_0}\sO_X((1-d)\infty)^n[t]\Big)[1]\simeq{\rm Lie}(E)(R);\\
%&&H^1\Big(X,\,\sO_X(-d\infty)^n[t]\x{t-\sum_{s=0}^rA_s\tau^s}\sO_X((1-d)\infty)^n[t]\Big)\simeq\exp_E^{-1}({\rm Lie}(E)(R));\\
%&&H^2\Big(X,\,\sO_X(-d\infty)^n[t]\x{t-\sum_{s=0}^rA_s\tau^s}\sO_X((1-d)\infty)^n[t]\Big)\simeq H^1(E/R);\\
%&&H^s\Big(X,\,\sO_X(-d\infty)^n[t]\x{t-\sum_{s=0}^rA_s\tau^s}\sO_X((1-d)\infty)^n[t]\Big)=0\;\hbox{ for } s\neq 1,\,2.
%\end{eqnarray*}
\bibliographystyle{plain}

%\bibliographystyle{amsplain}
%\bibliography{}
\end{document}